\documentclass[10pt,a4paper]{article}

\usepackage{amssymb}
\newtheorem{thm}{Theorem}
\newtheorem{defn}{Definition}
\newtheorem{cor}{Corollary}
\newtheorem{lem}{Lemma}
\usepackage{color}
\usepackage{colordvi}
\newcommand{\redtext}[1]{{\tt\textcolor{red}{#1}}}
\newcommand{\bluetext}[1]{{\tt\textcolor{red}{#1}}}

\newtheorem{prop}{Proposition}
\newtheorem{ex}{Example}
\newtheorem{rk}{Remark}
\usepackage{amssymb}

\def \endsquare{ $\sqcup \!\!\!\! \sqcap$ }

\def \M {{\! \rm \ I \!\! M}}

\def \E {{\! \rm \ I \!\! E}}
\def \B {{\! \rm \ I \!\! B}}
\def \A {{\! \rm \ I \!\! A}}

\def\otherterm#1{{\it#1}}

\def \Ci {{C^\infty}}
\def \l {{\lambda}}
\def \Cl {{C\ell}}

\def\res{{\rm res}}

\newcommand{\pdo}{$\Psi{\rm do}-$}
\newcommand{\pdos}{$\Psi{\rm dos}-$}

\newcommand{\tr}{{\rm tr}}

\newcommand{\AdP}{{\rm Ad}\ P}
\newcommand{\R}{\mathbb R}
\newcommand{\Z}{\mathbb Z}
\newcommand{\C}{\mathbb C}
\newcommand{\N}{\mathbb N}
\newcommand{\e}{\varepsilon}

\newcommand{\Ad}{{\rm Ad}}

\begin{document}

\title{Renormalised Chern-Weil forms  associated with families of
Dirac operators}

\author{ Jouko MICKELSSON, Sylvie
PAYCHA }

\maketitle
\section*{Abstract } 
We provide local expressions for  Chern-Weil type forms built from superconnections associated
with families of Dirac operators previously investigated in \cite{Sc} and later
in \cite{PS1}. \\
When the underlying fibration of manifolds
is trivial,  the even degree forms can be interpreted as renormalised Chern-Weil forms
in as far as they coincide with regularised Chern-Weil forms up to residue correction
terms.  Similarly, a new formula for the curvature of the local fermionic
vacuum line bundles  is derived using  a residue correction
term added to the naive curvature formula. \\
 We interpret the odd degree Chern-Weil type forms built from superconnections
 as Wodzicki residues and establish a
transgression formula along the lines of known transgression formulae for
$\eta$-forms. \\

\section*{Introduction}
Chern-Weil formalism in finite dimensions assigns to a connection $\nabla $ on
a principal bundle $P\to B$ over a manifold $B$ a form $f(\nabla)$ on $B$ with values in  the
adjoint bundle Ad $P$:
\begin{eqnarray*}
f: {\cal C}(P)&\to & \Omega(B, \Ad P)\\
 \nabla &\mapsto& f(\nabla),\\
\end{eqnarray*}
which is closed with de Rham cohomology class independent of the choice of
connection. Here $
{\cal C}(P)$ is the space of connections on $P$ and $\Omega(B, W)$ the space
of differential forms with values in a vector bundle $W$ over $B$.\\ \\
In the context of  \pdo bundles- i.e.
bundles with structure group the group of zero order invertible
pseudodifferential operators  $\Cl_0^*$- the trace used in finite dimensions to build maps $f_j(\nabla)= {\rm
  tr}(\nabla^{2j})$ can be replaced  by  two \footnote{the only two up
  to linear combinations \cite{LP}} natural traces on the algebra $\Cl_0$ of
zero order pseudodifferential operators, namely the Wodzikci residue and the
leading symbol trace. Such constructions were investigated in \cite{PR}  and lead to
maps which project down to quotient connections $\bar \nabla$ on the quotient
bundle $\bar P$ with structure group $\Cl_0^*/ \left(1+ \Cl_{-\infty}\right)^*$ where
$\Cl_{-\infty}$ is the algebra of smoothing operators and $\left(1+
  \Cl_{-\infty}\right)^*$ the group of invertibles. In other words, they
project down to   maps:
\begin{eqnarray*}
\bar f: {\cal C}(\bar P)&\to & \Omega(B, \Ad\bar  P)\\
 \nabla &\mapsto& \bar f(\bar \nabla).\\
\end{eqnarray*}
We call such maps {\it local} in as far as they are insensitive ``to smoothing perturbations''. 
\\ \\
In contrast, on a  principal bundle with structure group $\left(1+
  \Cl_{-\infty}\right)^*\subset \Cl_0^*$ one can mimic the
ordinary Chern-Weil construction to build Chern classes using the ordinary
trace on $\Cl_{-\infty}$.  We are concerned in this paper with possible
extensions of these Chern forms to  \pdo bundles.  Since the ordinary trace on
$\Cl_{-\infty}$ extends to linear functionals on $\Cl_0$ obtained from
regularised (or weighted) traces, one might want to try to extend the ordinary
Chern-Weil constructions to \pdo bundles using these regularised traces.  Such
issues were addressed in \cite{PR}; the fact that regularised traces  do not yield genuine traces gives rise to obstructions to 
carrying out the Chern-Weil construction since  the  regularised Chern forms
obtained from regularised traces are not closed.  However, it is useful to
keep in mind that  the obstruction to their closedness can 
be expressed in terms of  local maps  in the above sense. \\ \\
In this paper, we discuss ways to ``renormalise'' the regularised Chern forms
by adding  to them local maps in order to turn them into closed forms with de
Rham classes  independent of the connection. To do so, we compare
them with Chern forms  previously
investigated in \cite{Sc} and later \cite{PS1}, which are  built from
superconnections; in some cases they differ by local expressions so that  a 
renormalisation procedure  can indeed be carried out  adding local counterterms. More
precisely, letting (say  in the $\Z_2$-graded case)  $\A= D+ \nabla$ be a
superconnection associated with a Dirac operator $D$,  then the expression
$${\rm tr}^{D^2} (\nabla^{2j})- {\rm tr}^{\A^2}(\A^{2j})_{[2j]}$$ 
-which compares the naive infinite dimensional analog  ${\rm tr}^{D^2}
(\nabla^{2j})$ of the finite dimensional Chern form ${\rm tr} (\nabla^{2j})$
and the 
 closed form ${\rm tr}^{\A^2}(\A^{2j})_{[2j]}$ built from the super connection-  is local in the above sense. Here ${\rm
  tr}^{D^2}(B):= {\rm fp}_{z=0} {\rm TR} (B (D^2+\pi_D)^{-z})$ is the
$D^2$-weighted (or $\zeta$-regularised) trace of $B$ obtained as the finite
part at $z=0$ of the 
meromorphic expansion ${\rm TR} (B (D^2+\pi_D)^{-z})$ where $B$ is a form-valued 
pseudodifferential operator and TR the canonical trace on non integer order
pseudodifferential operators \cite{KV}.
$\pi_D$  stands for the orthogonal projection onto
the kernel of $D.$
\\
  This ``renormalisation''  procedure applies to the geometric setup
  corresponding to families of Dirac operators associated with a trivial
  fibration of manifolds (see Theorem \ref{thm:renChern}). \\
\\ In the case of  a family of Dirac operators
  associated with a general fibration of manifolds, such a straightforward
  ``renormalisation procedure'' is not possible  due to the presence of an
  extra curvature term arising from a  horizontal distribution 
on the fibration. Indeed, 
  the Chern-Weil forms associated
  with a superconnection  then differs from  a weighted Chern form by (a
  priori) non local
  terms  involving this extra curvature term.   \\For   a family
  of Dirac operators associated with a general fibration of spin manifolds
  $\pi: \M\to B$,
  on the grounds of the family index theorem, we identify   Chern forms
  associated with the superconnection with  form  components of $\int_{\M/B}
  \hat A (\M /B) \wedge {\rm ch} (\E_{\M/B})$ where $\E\to \M$ is a vector
  bundle over $\M$.
The $j$-th Chern form associated with a superconnection $\A$ introduced in
\cite{PS1} (following ideas of \cite{Sc}) has $2j$-form part  (see
  Theorem \ref{thm:localChernforms}): 
$$ {\rm str}^{\A^2} \left( \A^{2j}\right)_{[2j]}= \frac{(-1)^j j!}{(2i\pi)^{\frac{n}{2}}} \, \left( \int_{\M/B} \hat A(\M/B) \wedge
  {\rm ch}(\E_{\M/B}) \right)_{[2j]}.$$
On the grounds of the previous discussion, when the fibration is trivial,
it differs from renormalised weighted Chern forms by local terms. As it could be
expected in analogy with the finite dimensional situation, in the graded case,
the first
Chern form ${\rm str}^{\A^2}(\A^{2})_{[2]}$ turns out to be  proportional to
the curvature of the determinant bundle associated with the family of Dirac
operators from which the superconnection is built. \\
But there is also an  $j$-th  {\it residue} Chern form associated with a
superconnection $\A$ (which is new to our knowledge) the $2j-1$-th form part
of which which reads (see Theorem \ref{thm:localChernforms}):
$$
 {\rm sres}  \,  \left(\vert \A\vert^{2j-1}\right)_{[2j-1]}= \sqrt \pi  \frac{(-1)^j (2j-1)!!}{(2i\pi)^{\frac{n+1}{2}}\, 2^{j-1}} \, \left(\int_{\M/B} \hat A(\M/B) \wedge
  {\rm ch}(\E_{\M/B}))\right)_{[2j-1]}.$$
In the non graded  case, the second {\it residue} Chern form  (i.e. for $j=2$) turns out to be  proportional to the curvature
of the gerbe  associated with the family of Dirac
operators from which the superconnection is built, which was investigated by
Lott \cite{L}. \\  Following a similar scheme to that of Lott
  \footnote{We  derive a complete proof clarifying some
    steps in  Lott's proof.
  Our proof is carried out for operators $D(\lambda)=D-\lambda \, I$ (which are differential
    operators) but it easily extands to operators  $D_\alpha=
    D+h_\alpha(D)$ (which are pseudo-differential operators)  used by Lott where $h_\alpha$
    is a smooth function with compact support.}
 we derive a
transgression formula for the $j$-th residue Chern form (see Theorem \ref{thm:transgression})
$$
 {\rm sres} \left( \vert
 \A_\lambda\vert^{2j-1}\right)_{[2j-1 ]}=a_j \cdot  
  d\, \left(\tilde\eta_\lambda\right)_{[2j-2]},$$
using the $\eta$-invariant $\tilde\eta_\l$  (see \cite{BC},\cite{L})  associated with  a family of invertible Dirac type operators $D(\lambda)=
D-\lambda I$. These perturbed operators  differ from that of Lott but match physicists'
needs.\\ \\ 
The relation to gauge anomalies is explained in the last  section of the
paper. In particular, a new formula for the curvature of the local fermionic
vacuum line bundles (see Theorem \ref{thm:linebundlecurv}) is derived using  a residue correction
term added to the naive curvature formula (see Theorem \ref{thm:renormgauge}), the latter coming  by analogy  from the geometry of finite-dimensional Grassmann manifolds,
replacing the finite-dimensional trace by a weighted trace.  \\

\section{The geometric setup} 
Let  $\pi: E\to M$ be a vector bundle over a closed manifold $M$. $\Cl_0(M, E)$ denotes the Fr\'echet Lie algebra of $0$-order classical pseudo-differential
 operators (\pdo s) acting on smooth sections of $E$ and $\Cl_0^*(M, E)$ the Fr\'echet Lie group of invertible $0$-order classical pseudo-differential operators. \\ 
Let $P\to B$ be a $G=\Cl_0^*(M, E)$ principal bundle and $\Ad P= P\times_G \Cl_0(M,E)$ the adjoint
bundle, so that  locally, $\Ad P_{\vert_U}\simeq U\times \Cl_0(M, E)$. We equip $P$ with a connection $1$-form $\Theta: TP\to \Cl_0(M, E)$ which induces
 a connection $\nabla^\Ad$ on $\Ad P$. In   local coordinates we have $\nabla^\Ad= d+ [\Theta, \cdot]$ with $\Theta$ the above $\Cl_0(M, E)$-valued one form.\\
 \\ 
 The dual bundle $\Ad P^\star$  to $\Ad P$ comes equipped with the dual connection $\left(\nabla^\Ad\right)^\star$ defined for any section $\lambda$ of $\AdP^\star$
 and any section $\sigma$ of $\Ad P$ by $$d\lambda(\sigma)= \left(\left(\nabla^\Ad\right)^\star\lambda\right)(\sigma)+ \lambda(\nabla^\Ad \sigma).$$
On the other hand, $G$ acts on the space $C^\infty(M,E)$ of smooth sections of $E$  and  the associated vector bundle ${\cal E}= P\times_G C^\infty(M,E)$  comes
 equipped with the connection $\nabla$, locally of the form $\nabla= d+\Theta$. Then, locally $\nabla^\Ad= d+[\Theta, 
\cdot]$ and $\left(\nabla^\Ad\right)^*= d-[\Theta, \cdot].$ It is therefore convenient to write $\nabla^\Ad\sigma=[\nabla, \sigma]$ for any section 
$\sigma$ of $\Ad P$  and  $\left(\nabla^\Ad\right)^*\lambda=[\nabla, \lambda]$ for any section $\lambda$ of $\Ad P^\star$. With these notations we have:
$$d \left(\lambda(\sigma)\right)= [\nabla, \lambda] (\sigma)+ \lambda([\nabla ,\sigma]).$$
\\ \\
 The group  $\left(1+
  \Cl_{-\infty}(M,E)\right)^*$, where $\Cl_{-\infty}(M, E)$ denotes the algebra of smoothing operators,    is a normal subgroup of  $\Cl_0^*(M, E)$. Quotienting 
$\Cl_0^*(M, E)$  by  $\left(1+\Cl_{-\infty}(M, E)\right)^*$ yields   quotient  bundles
$\bar P\to B$ and $\bar {\cal E}=  \bar P\times_{\bar G} C^\infty(M,E)$ with structure group  $$\bar G:= \Cl_0^*(M, E)/\left(1+\Cl_{-\infty}(M, E)\right)^*$$
equipped with the induced connection $\bar \nabla$.\\
Let ${\cal C}(P)$ and ${\cal C}(\bar P)$ denote the space of connections on $P$ and $\bar P$. 
\begin{defn} We call a map \begin{eqnarray*}
f: {\cal C}(P)&\to &\Omega(B,{\Cl }( P))\\
\nabla &\mapsto& f(\nabla)
\end{eqnarray*}
{\rm local}  whenever it projects down to:
\begin{eqnarray*}
\bar f: {\cal C}(\bar P)&\to &\Omega(B,{\Cl }(\bar  P))\\
\nabla &\mapsto& \bar f(\bar \nabla).
\end{eqnarray*}
\end{defn}
Let  $\Cl ( {\cal E})= P\times_G \Cl(M, E)$ denote the bundle of classical pseudo-differential operators with fibre the whole algebra $\Cl(M,E)$ of classical
 pseudo-differential operators acting on sections of $E$. Clearly, $\Ad P\subset \Cl({\cal E})$ is a subbundle of $\Cl({\cal E})$.
\\ \\
In view of the following constructions, it is useful to mention that when $M=\{*\}$ is a point, then $E=V$ is a vector space,  $\Cl(M, E)= \Cl_0(M, E)= {\rm End}(V)$,  
$\Cl^*(M, E)= \Cl^*_0(M, E)= {\rm GL}(V)$ so that $P\to B$ boils down to an ordinary $GL(V)$-principal bundle and both $\Cl({\cal E})\to B$ and $\Ad P\to B$ boil down to
 its adjoint bundle $\Ad P= P\times_G {\rm End}(V)$ for the adjoint  action of $GL(V)$ on ${\rm End}(V)$. Thus, \pdo bundles can be seen as natural generalisations of 
ordinary principal bundles.

\section{$Q$-weighted traces (a short review)} 
A first atempt to generalise to \pdo bundles the construction of Chern-Weil
forms on ordinary bundles, is to use regularised   (or weighted) traces of
powers of the curvature as an Ersatz for ordinary traces of powers of the
curvature  which provide representatives of Chern-Weil classes in finite
dimensions \cite{PR}. We give a brief review of weighted traces of classical
pseudo-differential operators. \\ 
Let $Q\in \Cl(M, E)$ be an invertible admissible  elliptic operator of positive order $q$, where by admissible
we mean that its leading symbol admits a spectral cut $\theta$ \footnote{An operator $Q\in \Cl(M, E)$ of positive order  is called {\it
admissible } if there is a proper subsector of $\C$ with vertex 0
which contains the spectrum of the leading symbol $\sigma_L(Q)$ of
$Q$. Then there is a half line $L_\theta=\{r\, e^{i\theta}, r>0\}$ (a
spectral cut) with vertex $0$ and  determined by an Agmon angle
$\theta$ which does not intersect the spectrum of $Q$. If $Q$ is invertible,
then $\bar L_\theta=\{r\, e^{i\theta}, r\geq 0\}$ does not intersect the
spectrum of $Q$.}. If $Q$ is not invertible, we replace it by $Q+\pi_Q$
where  $\pi_Q$ is the orthogonal projection onto the kernel of $Q$. \\
An invertible admissible elliptic operator $Q$ has  complex powers $$Q_\theta^{z}=\frac{1}{2i\pi} \int_{\Gamma_\theta}
\lambda^{z} (Q-\lambda)^{-1} \, d\lambda$$
where $\Gamma_\theta$ is a contour around the spectral cut 
 and hence its
logarithm $\log_\theta Q=\frac{d}{dz}Q^{z}\vert_{z=0}$ which is not classical
anymore. In applications to follow, $Q$ is non negative self-adjoint so that
$\theta=\pi$ can be chosen as a spectral cut.  We shall henceforth drop
out the explicit mention of the spectral cut writing simply $Q^{-z}$ and $\log
Q$.  Recall that for any $A\in\Cl(M,E)$ and provided $Q$ has positive
  order, the map $z\mapsto {\rm TR} \left( A \,
  Q^{-z}\right)$ is meromorphic with simple pole at $0$ and its finite
part at $0$ $${\rm tr}^Q(A):= {\rm fp}_{z\to 0} {\rm TR} \left( A \,
  Q^{-z}\right)$$ is called the $Q$-weighted (or $\zeta$-regularised) trace of $A$. Here  TR is the canonical trace on non integer order classical \pdo
s [KV]. Even though it is not cyclic and hence not a genuine trace,  the
$Q$-weighted trace deserves the name of a trace in as far as it coincides with
the ordinary trace on trace-class operators and hence on  $\Cl_{-\infty}(M, E)$ and therefore extends it to  a
linear  map on  $\Cl(M,E)$.  \\
In contrast,  the Wodzicki residue  defined for $A\in \Cl(M, E)$ by
$${\rm res}(A)=\frac{1}{(2\pi)^n} \int_{S^*M} {\rm tr}_x \left( \sigma_A(x,
  \xi)\right)_{-n} \, dx \, d_S\xi$$ vanishes on trace-class operators and hence on
 $\Cl_{-\infty}(M, E)$. 
Here $S^*M$ stands for the cotangent unit sphere, $d_S\xi$ the canonical
volume measure on $S M$,  $\sigma_A$ is the symbol of $A$, ${\rm tr}_x$ the fibrewise trace and  the subscript $-n$ stands for the $-n$ (positively) homogeneous part of the 
symbol. \\ 
When $A$ is a differential operator we have \cite{PS2}:
\begin{equation}\label{eq:trQres}
{\rm tr}^Q(A)= -\frac{1}{q} {\rm res} \left(A\, \log Q\right)\end{equation}
where the residue on the r.h.s is defined by the above formula in spite of
$A\log Q$  not being  classical anymore.  The fact that $A$ is differential
ensures that the residue is well-defined. \\
In general ${\rm tr}^Q(A)$ depends on $Q$ for a given   $A\in \Cl(M,
E)$; given  two weights $Q_1, Q_2\in \Cl(M, E)$ with positive orders $q_1, q_2$ and same
spectral cut we have:
\begin{equation}\label{eq:trQ1Q2}
{\rm tr}^{Q_1}(A)-{\rm tr}^{Q_2}(A)= - {\rm res} \left(A\,\left( \frac{\log
      Q_1}{q_1}- \frac{\log
      Q_2}{q_2}\right)\right).\end{equation}
Also,  ${\rm tr}^Q$ is not cyclic: the obstruction to the cyclicity of ${\rm tr}^Q$ is measured by a Wodzicki
residue: 
\begin{equation}\label{eq:coboundary}
{\rm tr}^Q([A,B])= -\frac{1}{q} {\rm res} \left(A\,[B, \log
  Q]\right),\end{equation}
where now the residue is applied to a genuine classical operator since the
bracket $[B, \log
  Q]$ is classical. \\

We shall need the following technical lemma (see \cite{BGV} Lemma 9.35).
\begin{lem}\label{lem:Mellin} Let  $f$ be a smooth function on $]0, +\infty[$
  with  asymptotic behaviour at 0  of the type
$$f(\e)\sim_{\e\to 0} \sum_{j=0}^\infty a_j \e^{\alpha-j}$$ for some  real number
$\alpha$ (depending on $f$) and such that for large enough $\e$, $$\vert f(\e)\vert \leq C
e^{-\e \, \lambda}$$ for some $\lambda>0$, $C>0$. Then   its  Mellin transform  
 $$z\mapsto {\cal M}(f)(z):= \frac{1}{\Gamma(z)} \int_0^\infty \e^{z-1} f(\e) \,
 dt$$ defines  a meromorphic map on the complex plane (which  turns ut to be  holomorphic at $z=0$)  and  $${\rm fp}_{\e=0} f(\e)= {\rm fp}_{z=0} 
{\cal M}(f)(z)={\cal M}(f)(0).$$In particular,
 if $f(\e)= \sqrt \e \, g(\e)$ then 
$$ {\rm fp}_{\e=0} f(\e)= \sqrt \pi\, {\rm res}_{z=0} \left( {\cal M}(g)(z+\frac{1}{2})\right).$$
\end{lem}
{\bf Proof:}
The first part of the lemma is well known (see e.g. \cite{BGV} Lemma 9.35). Let us check the formula relating finite parts of $f(\e)=\sqrt{\e} g(\e)$ and its
Mellin transform.
\begin{eqnarray*}
 {\rm fp}_{\e=0} f(\e)
&=& {\rm fp}_{z=0} {\cal M}(f)(z)\\
&=& {\rm fp}_{z=0} \frac{1}{\Gamma(z)} \int_0^\infty \sqrt \e\,  \e^{z-1} g(\e)d\e\\
&=&  {\rm fp}_{z=0} \frac{1}{\Gamma(z)} \int_0^\infty  \e^{z+\frac{1}{2}-1} g(\e)d\e\\
&=&  {\rm fp}_{z=0} \frac{\Gamma(z+\frac{1}{2})}{\Gamma(z)}{\cal M}(g)(z+\frac{1}{2})\\
&=&  \Gamma(\frac{1}{2})\,{\rm fp}_{z=0}\left( z\, {\cal M}(g)(z+\frac{1}{2})\right)\\
&=& \sqrt \pi\,  {\rm res}_{z=0} \left( {\cal M}(g)(z+\frac{1}{2})\right).\\
\end{eqnarray*} 
 \endsquare
\\ \\ The Mellin transform provides a stepping stone  between heat-kernel regularisation and $\zeta$-regularisation methods: 
\begin{prop}\label{prop:trQMellin}
For any  $Q\in \Cl(M, E)$ non negative  self-adjoint elliptic and any $A\in \Cl(M, E)$ with  vanishing 
Wodzicki residue:$${\rm
  tr}^Q(A)= {\rm fp}_{\e = 0} {\rm tr}\left(A \, e^{-\e Q}\right).$$
\end{prop}
{\bf Proof:} This follows from  Lemma \ref{lem:Mellin} applied  to $f(\e)= {\rm tr} (A\, e^{-\e Q})$.\endsquare  
\section{$Q$-weighted Chern forms}
 We define  weighted
 Chern-Weil forms as in \cite{PR} and  briefly  recall  the obstructions
 to the closedness.\\
 Weighted traces  extend to \pdo valued forms in the following manner. 
Given a \pdo vector bundle  ${\cal E}$,   $Q$ is   a section of  $\Cl ( {\cal E})$ which is  elliptic,
admissible  and has  positive constant order $q$. Note that  these properties,
ellipticity, admissibility and  of  constant order $q$ are invariant under the
adjoint action of the group $\Cl^*(M, E)$ of invertible classical
pseudo-differential operators. The definition of the $Q$-weighted trace and
the Wodzicki residue then extend to \pdo valued forms setting for $b\in
U\subset B$ and $\alpha\otimes  A  \in\Omega\left(U,\Cl( {\cal E})\right)$, with $\alpha \in \Omega(U)$, $A\in
\Ci(U, \Cl(  {\cal E}))$:
$${\rm tr}^Q(\alpha\otimes A)(b):= \alpha(b) \otimes{\rm tr}^{Q_b}( A(b));
\quad {\rm res}(\alpha\otimes  A)(b):= \alpha(b) \otimes {\rm res}(
A(b)).$$ Properties (\ref{eq:trQres}), (\ref{eq:trQ1Q2})
(\ref{eq:coboundary}) extend in a straighforward manner to forms:
\begin{eqnarray}\label{eq:trQforms}
 {\rm tr}^Q(\alpha)&=& -\frac{1}{q} \, {\rm res}(\alpha \, \log Q)\nonumber\\
{\rm tr}^{Q_1}(\alpha)- {\rm tr}^{Q_2}(\alpha)&=& - \, {\rm res}\left(\alpha \left(
\,\frac{ \log Q_1}{q_1}- \frac{ \log Q_2}{q_2}\right)\right)\nonumber\\
 {\rm tr}^Q([\alpha,\beta])&=& -\frac{1}{q} \, {\rm res}(\alpha \, [\beta,\log
 Q]), 
\end{eqnarray}
where the first identity holds whenever $\alpha$ is a differential operator
valued form whereas the others hold for any $\Cl\left({\cal E}\right)$-valued forms
$\alpha, \beta$.\\
 The
Wodzicki residue commutes with differentiation 
$[\nabla, {\rm res}]=0$ whereas weighted traces do not. The obstruction is
measured in terms of a Wodzicki residue. Indeed, it follows from 
(\ref{eq:trQforms}) that locally,  $[d\, {\rm tr}^Q](\alpha)=- \frac{1}{q} {\rm res}
( \alpha\, d\log Q)$ as a result of which, writing $\nabla^{\rm Ad}=
d+[\theta, \cdot]$ in local coordinates, we have:
\begin{eqnarray} \label{eq:nablatrQ}
 [\nabla,  {\rm tr}^
Q](\alpha)&=& d\left({\rm tr}^
Q(\alpha)\right)-{\rm tr}^
Q([\nabla, \alpha])\nonumber \\
&=& d\left({\rm tr}^
Q(\alpha)\right)-{\rm tr}^
Q(d\, \alpha)-{\rm tr}^
Q([\theta,  \alpha])\nonumber\\
&=& [ d\,  {\rm tr}^
Q]( \alpha)- \frac{1}{q} \, {\rm res}\left(   \alpha\, [\theta,\log
  Q]\right)\nonumber\\
&=& - \frac{1}{q} {\rm res}
( \alpha\, d\log Q)- \frac{1}{q} \, {\rm res}\left(   \alpha\, [\theta,\log
  Q]\right)\nonumber\\
&=&  -\frac{1}{q} {\rm res}\left( \alpha\, [\nabla, \log Q]\right)\quad \forall \alpha\in \Omega\left(B, \Cl({\cal E})\right).
\end{eqnarray}

The curvature   $\Omega =\left(\nabla^\Ad\right)^2$ 
 of $\nabla^\Ad$ lies in $\Omega^{2}\left(P,\Cl\left( {\cal E}\right)\right)$ so
 that   the $Q$-weighted trace
${\rm tr}^Q(\Omega^i)$
 defines a $2i$-form on $B$.  
  The following proposition tells us that the obstruction to the closedness is local in the sense of the above definition.
\begin{prop}The exterior differential of the weighted Chern-Weil form ${\rm tr}^Q(\Omega^i)$ is local, i.e. of the form $d{\rm tr}^Q(\Omega^i)=f_i(\bar \nabla)$ for some $\bar f_i: {\cal C}(\bar P)\to \Omega(\Cl(\bar P))$.
\end{prop} 
{\bf Proof:}
Since  $\nabla^{\Ad}(\Omega)=[\nabla, \Omega^ i]=0$, by (\ref{eq:nablatrQ}) we have
\begin{eqnarray*}
d {\rm tr}^Q(\Omega^i)&=& [\nabla,  {\rm tr}^
Q](\Omega^i)+ {\rm tr}^
Q([\nabla, \Omega^i])\\
&=& [\nabla,  {\rm tr}^
Q](\Omega^i)+\sum_{j=1}^i {\rm tr}^
Q(\Omega^j\,[\nabla,\Omega]\,\Omega^{i-j})\\
& =& [\nabla, {\rm tr}^
Q](\Omega^i)\\
&=& -\frac{1}{q} {\rm
  res}\left(\Omega^i\,  [\nabla,\log  Q]\right).
\end{eqnarray*}
Since the Wodzicki residue vanishes on smoothing operators,
$$d\, {\rm tr}^Q(\Omega^i)=d\,  {\rm tr}^Q(\bar \Omega^i)= \bar f(\bar \nabla),$$ and hence the locality property of the obstruction to the closedness.
 \endsquare
.

\section{From  superconnections to Chern-Weil type  forms}
We review and extend constructions of Chern-Weil type forms carried out in
\cite{PS1} using superconnections.  Let ${\cal E}$ be a vector bundle associated with a
\pdo principal bundle $P$ as before. 
\subsection{Chern forms associated with superconnections}
\begin{itemize}
\item {\bf The $\Z_2$-graded case:}
Let us  assume that ${\cal E}= {\cal E}^+ \oplus {\cal E}^-$ is  a
$\Z_2$-graded super bundle on $B$. The canonical trace TR for non integer
order operators  in $\Cl({\cal
  E})$ is replaced by the super canonical trace sTR whereas weighted traces ${\rm tr}^Q$ for operators in $\Cl({\cal
  E})$ are  replaced by weighted supertraces ${\rm str}^Q$ with respect to even weights
 $Q= Q^+\oplus Q^-$. They vanish on odd \pdos and give the difference of
 weighted traces on even \pdos:
$$ {\rm sTR}(A):= {\rm
  TR}(A^{++})- {\rm
  TR}(A^{--}); \quad {\rm str}^Q(A):= {\rm
  tr}^{Q^+}(A^{++})- {\rm
  tr}^{Q^-}(A^{--})$$ with obvious notations.\\
This grading combined with the $\Z_2$-grading on forms $\Omega(B, {\cal E})= \Omega^{ev} (B, {\cal E})\oplus \Omega^{od} (B, {\cal E})$ gives rise to:
$$\Omega^+(B, {\cal E})=  \Omega^{ev} (B, {\cal E}^+)\oplus \Omega^{od} (B,
{\cal E}^-); \quad \Omega^-(B, {\cal E})=  \Omega^{od} (B, {\cal E}^+)\oplus
\Omega^{ev} (B, {\cal E}^-).$$
\begin{defn}\cite{BGV}
A superconnection is an odd-parity first order differential operator $$\A:
\Omega^{{}^+_-}(B, {\cal E})\to \Omega^{{}^-_+}(B, {\cal E})$$ which satisfies
the (graded) Leibniz  rule. If $\alpha \in \Omega(U)$ for some open subset
$U\subset B$ and $B\in
\Ci\left(U,\Cl( {\cal E})\right)$ then 
$$\A (\alpha \otimes B)= d\alpha \otimes B+ (-1)^{\vert \alpha \vert}\alpha
\otimes  \A\, B.$$ 
\end{defn}
 The curvature 
$\A^2$ of a superconnection $\A$ on ${\cal E}$
lies in $\Omega(B,\Cl({\cal E}))$.  Following Quillen [Q] we say a superconnection
$\A$ on ${\cal E}$ is associated with a smooth family of elliptic differential operators
$\{D_b, b\in B\}$  whenever $\A_{[0]}= D$.
\begin{ex}
   $\A= D+\nabla$ defines a particular super connection associated with
  $D$ with  curvature 
$$\A^2= D^2 + \nabla^\Ad D+ \Omega=Q+ [\nabla, D]+ \Omega$$
where we have set $Q=D^2$. Here $[\nabla, D]=\nabla D+ D\nabla $ is the anticommutator. 
\end{ex}
\begin{rk}In the following we systematically use  graded commutators of operator valued forms: 
anticommutator for odd-odd forms and usual commutator otherwise.
\end{rk}
\item {\bf The non graded case:}
Let us  assume that ${\cal E} $ is an ordinary vector bundle on $B$. Following Quillen, we introduce an extra grading $\sigma$ such that
$\sigma^2=1$ and build the right $\C\oplus  \C\, \sigma$-module:
$$\Cl_\sigma \left({\cal E} \right):= \Cl \left({\cal E} \right)\bigotimes
\left( \C\oplus  \C\, \sigma\right).$$
Odd degree \pdos lie in  $\Cl \left({\cal E} \right)\bigotimes
\left(  \C\, \sigma\right)$ whereas $\Cl \left({\cal E} \right)$ is identifued
with even degree \pdos.  
The ordinary  canonical trace TR extends to non integer order operators in  
$\Cl_\sigma \left({\cal
     E}\right)$ by ${\rm sTR}(\alpha +\sigma \beta)= {\rm TR}(\beta)$ so
 that weighted traces ${\rm tr}^Q$ are  replaced by  $${\rm str}^{Q}(\alpha + \sigma \beta):=   {\rm tr}^{Q}(
  \beta).$$
These definitions
extend to the space $\Omega_\sigma \left(B, {\cal E} \right)$  of $\Cl_\sigma \left({\cal E} \right)$-valued forms on $B$ in a
straightforward manner.
\begin{defn}\cite{BGV}
A superconnection is a first order differential operator $$\A:
\Omega_{\sigma}(B, {\cal E})\to \Omega_\sigma(B, {\cal E})$$ which  commutes with $\sigma$ and
satisfies
the (graded) Leibniz  rule.
\end{defn}
As in the even case, it is associated with a family $\{D_b, b\in B\}$ of
elliptic differential
operators whenever $\A_{[0]}=D$. 
\begin{ex}$\A:= \sigma \, D+ \nabla$ is a particular  superconnection associated with
$D$, the curvature of which reads \footnote{Here $D$
commutes with $\sigma$ whereas $\nabla$
anticommutes with $\sigma$.}
$$\A^2= D^2 + \nabla^\Ad(\sigma  D)+ \Omega=Q+ [\nabla, \sigma\, D]+ \Omega$$
with $Q=(\sigma\, D)^ 2= D^ 2$ as before. 
\end{ex}
\end{itemize}
Let us recall  from  \cite{Sc} (see also \cite{PS1}) that  weighted traces can
be extended to include weights $\A^2$ which are
\pdo valued forms  and   analogs of Chern-forms can be constructed, which  turn out to be
closed.   Writing $$\A^2= \A_{[0]}^2+ \A_{[1]}^2+ \A_{[2]}^2= D^2+ \A^2_{[>0]}$$
where the subscript $[j]$ stands for the $j$-th degree part, and $[>0]$ for
non no  zero degree part,   can be  useful to derive explicit expansions in
increasing form degree. 
For example, 
 \begin{eqnarray}\label{eq:resolv}
 (  \lambda-\A^2)^{-1}&=& \left( \lambda- D^2-
     \A^2_{[>0]}\right)^{-1}\nonumber \\
&=& \sum_{j=0}^K  \left( \lambda- D^2\right)^{-1}
     \A^2_{[>0]} \left( \lambda- D^2 \right)^{-1}\cdots  \A^2_{[>0]} \left( \lambda- D^2\right)^{-1}\nonumber\\
&+&S_K (D^2, \A^2_{[>0]}, \lambda), 
\end{eqnarray}
where $S_K $ has form degree $>K$ and where $\A^2_{[>0]} \left( \lambda-
  D^2)\right)^{-1}$ arises $j$ times in the $j$-th term of the sum. By
convention the $j=0$ term reduces to $\left( \lambda-
  D^2)\right)^{-1}$.  Hence, for any positive integer $K$
 $$
 \left((  \lambda-\A^2)^{-1}\right)_{[K]}=  \sum_{j=0}^K  \left( \lambda- D^2\right)^{-1}
     \A^2_{[>0]} \left( \lambda- D^2)\right)^{-1}\cdots  \A^2_{[>0]} \left(
     \lambda- D^2\right)^{-1}$$
has a finite expansion in increasing form degree. 
Also, whenever $\A_{[0]}=D$ is invertible, so is  $\A^2$ invertible and its modulus $\vert \A\vert
  :=\left(\A^2\right)^{\frac{1}{2}}$   can be defined using a contour
  integration (cfr Section 2):
$$\vert \A\vert= \frac{i}{2\pi} \int_\Gamma
\sqrt \lambda\, (\A^2-\lambda)^{-1}\, d\lambda$$
where $\Gamma$ is a contour around the spectrum of $D^2$ which is a subset of
$\R^+$. \\ In general, we set 
$$\vert \A\vert:=\sqrt{\A^2+\pi_\A}$$ where $\pi_\A$ is the orthogonal
projection onto the kernel of $\A_{[0]}^2$. This  defines a form provided Ker $\A_{[0]}^2$
has constant dimension. 
\begin{rk}
Note that for any $\alpha\in \Omega(B, \Cl\left({\cal E})\right)$,
$${\rm str}^{\A^2} (\alpha)=  {\rm fp}_{z=0}  {\rm sTR} \left( \alpha  \,
 \left( \A^2+\pi_\A\right)^{-z}\right)= {\rm fp}_{\e=0}  {\rm str} \left( \alpha
  \, e^{-\e \A^2}\right)$$ since $\A^2$ is a differential operator-valued form  and hence has
vanishing Wodzicki residue. 
\end{rk}
 The following proposition extends results of \cite{PS1}.
 \begin{prop}\label{prop:closedChernforms}Let  $P$ be a polynomial function. Forms  ${\rm str}^{\A^2}
(P(\A^2))$   and ${\rm sres}\left(P(\vert  \A)\vert \right)$ 
associated with a superconnection $\A$ are closed. Their de Rham class is
  independent of the choice of connection  one forms $\A_ {[1]}$. \\

\end{prop}
\begin{rk} The residue form  ${\rm sres}\left(P(\vert  \A\vert
    )\right)$ is in fact insensitive to the projection $\pi_\A$ which is a
    smoothing  operator and hence does not affect the Wodzicki residue. 
\end{rk}
{\bf Proof:}
We extend   the argument used in \cite{PS1} for the closedness
of  forms  ${\rm str}^{\A^2}
(\A^{2i})$ to any ${\rm str}^{\A^2}
(P(\A))$. 
Equations  (\ref{eq:nablatrQ}) and (\ref{eq:coboundary}) extend replacing
 the connection  $\nabla$ by the
superconnection $\A$ and the weight  $Q$ by  the \pdo valued form $\A^2$ of
order $2$ \cite{PS1} and  we have
\begin{eqnarray*} d \, {\rm str}^{\A^2}\left(P(\A^2)\right)
&=& [\A,  {\rm str}^{\A^2}]
\left(P(\A^2)\right)+  {\rm str}^{\A^2}\left([\A, P(\A^2)]\right)\\
&=& -\frac{1}{2} {\rm sres} \left( P(\A^{2}) [\A,  \log( \A^2+\pi_\A)]\right)\\
&=&-\frac{1}{2}\, \frac{d}{dt}_{t=0}\left( \frac{i}{2\pi} \int_{\Gamma}  \lambda^t \,
 {\rm sres} \left( P(\A^{2})[\A, (\A^2+\pi_A-\lambda)^{-1}]\right) \, d\lambda\right)\\
 &=& 0
\end{eqnarray*}
since $[\A, (\A^2+\pi_A-\lambda)^{-1}]$ is smoothing. \\ 
 
Similarly, 
$$ d\,  {\rm sres}\left(P(\vert \A\vert)\right)
=   {\rm sres}\left([\A, P(\vert \A\vert)]\right)
=0$$
since $[\A,  \log P(\vert \A\vert)]=0$.\\
 The forms are therefore closed.\\ 
 Let us check that their de Rham classes are independent of the
choice of superconnection. Let $\A_t$
be a smooth one parameter family of superconnections then for any monomial
$P(\A^2)= \A^{2i}$

\begin{eqnarray}\label{eq:variation1}
&{}& \frac{d}{dt}\left( {\rm
      str}^{\A_t^2}\left(\A_t^{2i}\right)\right) \nonumber\\
&=&  \frac{d}{dt}\left({\rm fp}_{\e=0} {\rm str}\left(\A_t^{2i}\, e^{-\e
      \A_t^2}\right)\right)\nonumber\\
&=&  \sum_{j=1}^i{\rm fp}_{\e=0} \left({\rm str}\left(\A_t^{2(i-j-1)}\,[\A_t,
    \dot \A_t]\,\A_t^{2j}\,  e^{-\e
      \A_t^2}\right)-\e \,{\rm str}\left(\A_t^{2i}\,[\A_t,
    \dot \A_t]\, e^{-\e
      \A_t^2}\right)\right)\nonumber\\
&=& i\,{\rm fp}_{\e=0} \left({\rm str}\left([\A_t, \A_t^{2(i-1)}
    \dot \A_t\, e^{-\e
      \A_t^2}]\right)- \e \,  {\rm str}\left([\A_t, \A_t^{2i}\,
    \dot \A_t\, e^{-\e
      \A_t^2}]\right)\right)\nonumber\\
&=& i\,{\rm fp}_{\e=0} \left(d\, {\rm str}\left(\A_t^{2(i-1)}
    \dot \A_t\, e^{-\e
      \A_t^2}\right)- \e \, d\, {\rm str}\left(\A_t^{2i}\,
    \dot \A_t\, e^{-\e
      \A_t^2}\right)\right)\nonumber\\
&=& d\, \left[ i\,{\rm fp}_{\e=0} \left( {\rm str}\left(\A_t^{2(i-1)}
    \dot \A_t\, e^{-\e
      \A_t^2}\right)- \e \,  {\rm str}\left(\A_t^{2i}\,
    \dot \A_t\, e^{-\e
      \A_t^2}\right)\right)\right]
\end{eqnarray}

is exact.
Here we have used the fact that $\frac{d}{dt} \A_t^2= \dot \A_t\, \A_t + \A_t
\, \dot \A_t= [\A_t, \dot \A_t]$, the graded commutator of $\A_t$ with the
form \pdo valued form $\dot \A_t$. 
It follows that the de Rham class of ${\rm str}^{\A^2}
\left(\A^{2i}\right)$ and hence of  ${\rm str}^{\A^2}
\left(P(\A^{2})\right)$ is independent of the choice of connection.\\ \\
Similarly, since 
 $\vert \A\vert^j= \frac{i}{2\pi} \int_\Gamma \lambda^{\frac{j}{2}}\, (\A^2+\pi_A-\lambda)^{-1}\, d\lambda$ and since:
 \begin{eqnarray*}
\frac{d}{dt}
  (\A_t^2+\pi_A-\lambda)^{-1}&=& -  (\A_t^2+\pi_\A-\lambda)^{-1}\, \frac{d}{dt}\A_t^2\, 
  (\A_t^2+\pi_\A-\lambda)^{-1}\\
&=& -  (\A^2+\pi_\A-\lambda)^{-1}\, [\A_t, \, \dot \A_t]
  (\A_t^2+\pi_\A-\lambda)^{-1}\\
&=&- \left[\A_t,\, (\A_t^2+\pi_\A-\lambda)^{-1} \,\dot \A_t\, 
  (\A_t^2+\pi_\A-\lambda)^{-1}\right],\\
\end{eqnarray*}
it follows that the variation

\begin{eqnarray}\label{eq:variation2}
&{}& \frac{d}{dt}\left( {\rm sres}\left(\vert\A_t\vert^{j}\right)\right)\nonumber\\
&=& {\rm sres}\left(\frac{i}{2\pi}\left[\int_\Gamma \lambda^{\frac{j}{2}}\,  \frac{d}{dt}
  (\A_t^2+\pi_A-\lambda)^{-1}\, d\lambda\right] \right)\nonumber\\
&=&- \frac{i}{2\pi}\, {\rm sres}\left(\int_\Gamma \lambda^{\frac{j}{2}}\,  
 \left[\A_t,  (\A_t^2+\pi_\A-\lambda)^{-1} \,\dot \A_t
  (\A_t^2+\pi_\A-\lambda)^{-1}\right]\,d\lambda\right)\nonumber\\
&=&- \frac{i}{2\pi}\, {\rm sres}\left(\int_\Gamma \lambda^{\frac{j}{2}}\,  
 \left[\A_t+\pi_\A,  (\A_t^2+\pi_\A-\lambda)^{-1} \,\dot \A_t
  (\A_t^2+\pi_\A-\lambda)^{-1}\right]\,d\lambda\right)\nonumber\\
&{}& {\rm since} \quad {\rm res} [\pi_\A, \cdot]=0\nonumber\\
&=&- \frac{i}{2\pi}\, {\rm sres}\left( \left[\A_t+\pi_\A,\int_\Gamma  \lambda^{\frac{j}{2}}\,  
  \,\dot \A_t
  (\A_t^2+\pi_\A-\lambda)^{-2}\right]\right)\nonumber\\
&{}& {\rm since}\quad [\A_t+\pi_\A, (\A_t^2+\pi_\A-\lambda)^{-1}]=0 \nonumber\\
&=&- \frac{i}{2\pi}\, {\rm sres}\left( \left[\A_t,\int_\Gamma  \lambda^{\frac{j}{2}}\,  
  \,\dot \A_t\, 
  (\A_t^2+\pi_\A-\lambda)^{-2}\right]\right)\nonumber\\
&{}& {\rm since} \quad {\rm res} [\pi_\A, \cdot]=0\nonumber\\
&=&\frac{j}{2} \frac{i}{2\pi}\, {\rm sres}\left( \left[\A_t,\int_\Gamma  \lambda^{\frac{j-2}{2}}\,  
  \,\dot \A_t\, 
  (\A_t^2+\pi_\A-\lambda)^{-1}\right]\right)\nonumber\\
&=&\frac{j}{2} \,{\rm sres}\left( \left[\A_t, 
  \,\dot \A_t\,
(\A_t^2+\pi_\A)^{\frac{j-2}{2}}\right]\right)\nonumber\\
&=& \frac{j}{2} \, d\,  \left({\rm sres}\left(  
  \,\dot \A_t\,
\vert \A_t\vert^{j-1}\right)\right)
\end{eqnarray}
is also  exact, which ends the proof of the proposition. 
\endsquare

\subsection{ Chern-forms associated with superconnections $\A= \nabla+ D$}
 We now specialise to the case $\A_{[2]}=0$ and consider a superconnection
 $\A= D+\nabla$ in the graded  setup and  $\A= \sigma\, D+\nabla$ in the
 ungraded setup.  The following theorem compares the closed Chern-forms ${\rm str}^{\A^2} (\A^{2i})$ with the  (non closed in general)
 weighted Chern forms.  
\begin{thm} \label{thm:renChern}In the $\Z_2$-graded set up and  provided the superconnection $\A= D+\nabla$, the (closed)  Chern-forms ${\rm str}^{\A^2} (\A^{2j})_{[2j]}$ 
differ from the (non closed) $Q$-weighted Chern forms ${\rm str}^Q(\Omega^j)$  by a local map  i.e. 
$${\rm str}^{\A^2}\left(\A^{2j}\right)_{[2j]}- {\rm str}^Q(\Omega^j)=\bar f_j(\bar \nabla). $$
 for some $\bar f_j: {\cal C}(\bar P)\to \Omega(B, \bar P).$ \\
In the ungraded setup and provided the superconnection $\A= \sigma\,
D+\nabla$,  (closed)  Chern-forms ${\rm str}^{\A^2} (\A^{2j})_{[2j-1]}$ 
differ from the (non closed) $Q$-weighted forms ${\rm str}^Q(\Omega^{j-1}\, [\nabla, \sigma\, D])$  by a local map  i.e. 
$${\rm str}^{\A^2}\left(\A^{2j}\right)_{[2j-1]}-j\, {\rm str}^Q(\Omega^{j-1}\,
[\nabla, \sigma\, D])=\bar g_j(\bar \nabla)$$ 
for some $\bar g_j: {\cal C}(\bar P)\to \Omega(B, \bar P).$
\end{thm}
\begin{rk} This does not hold anymore if $\A_{[2]}\neq 0$ as can easily be
  seen from the proof below. When $\A_{[2]}= 0$, on the grounds of this
  proposition,   ${\rm
    str}^{\A^2}\left(\A^{2j}\right)_{[2j]}$ can be interpreted  as a renormalised version of $ {\rm str}^Q(\Omega^j)$.
\end{rk}
{\bf Proof:}
Let us observe  in the graded case that since $\A^{2}= Q+ [\nabla, D]+ \Omega$, we have:
$${\rm str}^Q\left(\A^{2j}\right)_{[2j]}= {\rm str}^Q\left(\Omega^{j}\right),$$
 and similarly in the ungraded case, we have   $${\rm str}^Q\left(\A^{2j}\right)_{[2j-1]}=j {\rm str}^Q\left(\Omega^{j-1}\, [\nabla,\sigma\, D]\right).$$  Using a
Campbell-Hausdorff formula for pseudo-differential operators \cite{O} combined with
 formula (\ref{eq:trQ1Q2}) extended to form valued weights, we have 
\begin{eqnarray*}
&{}& {\rm str}^{\A^2}\left(\A^{2j}\right)_{[2j]}\\
&=&  {\rm str}^Q\left(\A^{2j}\right)_{[2j]} +\left({\rm str}^{\A^2}\left(\A^{2j}\right)_{[2j]}-
 {\rm str}^Q\left(\A^{2j}\right)_{[2j]}\right)\\
&=& {\rm str}^Q\left(\Omega^{j}\right) -\frac{1}{2} {\rm sres}\left(\A^{2j}(\log
 \A^{2}- \log  Q)\right)_{[2j]}\\
&=& {\rm str}^Q(\Omega^j)-\frac{1}{2} {\rm sres}\left(\A^{2j}
\left(\log (1+Q^{-1}([\nabla, D]+ \Omega)) +[\log Q, \log (1+Q^{-1}([\nabla, D]+ \Omega))]+...\right) \right)_{[2j]}\\
&=& {\rm str}^Q(\Omega^{j-1})+ f_i( \nabla),
\end{eqnarray*}
with $f_j(\nabla)$ the Wodzicki residue of a polynomial expression in  $D$, 
$D^{-1}$
and $\nabla$ of total form degree $2j$. As  a Wodzicki residue, it is insensitive to a perturbation of the connection by a smoothing operator so that 
 $ f_j( \nabla)=\bar f_j(\bar \nabla)$. This shows that $${\rm str}^{\A^2}\left(\A^{2i}\right)_{[2j]}- {\rm str}^Q(\Omega^j)=\bar f_j(\bar \nabla)$$ is local. 
A similar computation shows that $${\rm str}^{\A^2}\left(\A^{2j}\right)_{[2j-1]}-j\, {\rm str}^Q(\Omega^j\, [\nabla, \sigma\, D])=\bar g_j(\bar \nabla)$$ is also local.
\endsquare 
\subsection{ Residue
    Chern forms as Wodzicki residues}
In order to derive an explicit expression for the residue Chern forms in terms
of Wodzicki residues, we
 borrow the following
notations from \cite{CoM} and \cite{H}.  For    $A$ in  $ Cl(M,
E)$   of order $a$, a given $\Delta\in \Cl(M, E)$ and  any $j\in \N$ we set: $$A^{(j)}:= {\rm ad}_{\Delta}^j(A),
\quad {\rm where} \quad  {\rm ad}_{\Delta}(B)=[\Delta, B],$$  so that $ A^{(0)}=A, \quad
 A^{(j+1)}= {\rm ad}_{\Delta}(A^{(j)})=[\Delta, A^{(j)}].$
When $\Delta $ of  order $2$ has scalar leading symbol then  $A^{(j)}$ has order $a+j+1$. 
\begin{prop} Let  $\A$ be a superconnection
associated with an operator $D$, the square of which  has scalar
leading symbol.  For any positive  integer $K$ 
\begin{eqnarray*}
 {\rm sres} \left( \vert \A\vert^{2j-1}\right)_{[K]}
&=&  \sum_{l=0}^K   \sum_{ k_1\geq
  0}\cdots  \sum_{ k_l\geq
  0}\frac{\left(\frac{2j-1}{2}\right)\cdots\left(\frac{2j-1}{2}-\vert k\vert
    -l\right) }{( k_1+\cdots +k_l +l)! }c(k_1, \cdots, k_l)\cdot \\
&{}& \cdot {\rm sres} \left(
\left(\A^2_{[>0]}\right)^{(k_1)}\left(\A^2_{[>0]}\right)^{(k_2)}\cdots
\left(\A^2_{[>0]} \right)^{(k_l)}\, \left(D^2\right)^{ \frac{2j-1}{2} -\vert k\vert-l}
\right)_{[K]},
\end{eqnarray*}
where  we set
$c(k_1)=1$ for any positive integer $k$ and where,  for a multi index $k=(k_1,\cdots, k_l)$ 
 for  $j>1$ we set 
$$c(k_1, \cdots, k_l)= \frac{(k_1+\cdots +k_l+l)!}{k_1!\cdots
  k_j!(k_1+k_2+1)\cdots (k_1+\cdots +k_{l-1}+l)}.$$  \\
In particular,
$$ {\rm sres} \left( \vert \A\vert\right)_{[1]}\\
=   \sum_{ k\geq
 0}\frac{\left(\frac{1}{2}\right)\cdots\left(\frac{1}{2}- k
    -1\right) }{( k +1)! }\cdot {\rm sres} \left(
\left(\A^2_{[1]}\right)^{(k)}\, \left(D^2\right)^{ \frac{1}{2} - k-1}
\right),
$$
\end{prop}

\begin{rk} If $D$ is a differential operator then $\A^2$ is a differential
  operator and ${\rm sres} \left(\vert \A\vert^{2j}\right)= {\rm sres}
  \left(\left( \A^2\right)^{j}\right)$ which is why we only consider odd
  powers.
\end{rk}
\begin{rk} Since the operator order of $\A^2_{[>0]}$ is no larger than $1$,
 $\left(\A^2_{[1]}\right)^{(k)}$ has  order $\leq 1+k$  and 
$\left(\A^2_{[>0]}\right)^{(k_2)}\cdots
\left(\A^2_{[>0]} \right)^{(k_l)}\, \left(D^2\right)^{ \frac{2j-1}{2} -\vert
  k\vert-l}$ has order $ \leq 2j-1 -\vert
  k\vert-l$ which decreases as $\vert k\vert $ or $l$ increases. Thus the Wodzicki residue
  vanishes for large enough $\vert k\vert $ or $l$ and  the seemingly
  infinite series in the proposition  is in fact finite.
\end{rk}
{\bf Proof:}
We introduce  notations borrowed from \cite{H} and \cite{CoM}. Let $T\in \Cl(M,E)$ and $T_k, k\in \N$
be operators in $\Cl(M, E)$  with decreasing order in $k$. Then 
$$ T\simeq \sum_{k\geq 0} T_k
\Longleftrightarrow  \forall N \in \N , \exists K(N)
\quad T- \sum_{k= 0}^{K(N)} T_k\in \Cl^{-N}(M, E).$$
With these notations,  for  any non negative  integer $h$ we have \cite{H} (see the proof of Proposition
4.14)
$$  (\lambda-D^2)^{-h} A\simeq \sum_{k\geq 0}\frac{ (h+k-1)!}{(h-1)!k!}A^{(k)}
(\lambda-D^2
)^{-h-k}. $$
As a result, the  $j$-th term in (\ref{eq:resolv})  reads: 
\begin{eqnarray*}
&{}& (\lambda-D^2)^{-1} \, \A^2_{[>0]}
\cdots  (\lambda-D^2)^{-1} \, \A^2_{[>0]}\,(\lambda-D^2)^{-1}\\
&\simeq& \sum_{k_1\geq 0} \left(\A^2_{[>0]}\right)^{(k_1)}
(\lambda-Q)^{-2-k_1}\A^2_{[>0]}
\cdots  (\lambda-D^2)^{-1} \A^2_{[>0]}\,(\lambda-D^2)^{-1}\\
&\simeq& \sum_{k_1\geq 0} \left(\A^2_{[>0]}\right)^{(k_1)}\sum_{k_2\geq 0} \frac{(-1)^{k_2} (k_1+k_2)!}{k_1!k_2!}\left(\A^2_{[>0]}\right)^{(k_2)}
(\lambda-D^2)^{-3-k_1-k_2} \A^2_{[>0]}\\
&{}&
\cdots  (\lambda-D^2)^{-1} \A^2_{[>0]}\,(\lambda-D^2)^{-1}\\
&\simeq& \sum_{k_1\geq 0} \left(\A^2_{[>0]}\right)^{(k_1)}\sum_{k_2\geq 0}
\frac{ (k_1+k_2)!}{k_1!k_2!}\left(\A^2_{[>0]}\right)^{(k_2)}\sum_{k_3\geq 0}\,
\\ & \cdot& \frac{(-1)^{k_3} 
(k_1+k_2+1)!}{(k_1+k_2+1)!k_3!}\left(\A^2_{[>0]}\right)^{(k_3)}
(\lambda-D^2)^{-4-k_1-k_2-k_3}  \A^2_{[>0]}\cdots\\
&{}&
\cdots  (\lambda-D^2)^{-1}  \A^2_{[>0]}\,(\lambda-D^2)^{-1}\\
&\simeq& \sum_{\vert k\vert\geq 0}c(k_1, \cdots, k_j)
\left(\A^2_{[>0]}\right)^{(k_1)}\left(\A^2_{[>0]}\right)^{(k_2)}\cdots
\left(\A^2_{[>0]}\right)^{(k_j)}\,(\lambda-D^2)^{-\vert k\vert-j-1}.\\
\end{eqnarray*}
Letting  $\Gamma$  be  a contour around the spectrum ${\rm spec}(D^2)\subset
\R^+$, it follows that for any positive integer $K$:
\begin{eqnarray*}
 &{}& \left( \left(\A^2\right)^{\frac{2i-1}{2}}\right)_{[K]}\\
&=& \sum_{j=0}^K \frac{1}{2i\pi} \int_{\Gamma} \, \lambda^{\frac{2i-1}{2}} \, 
 \left((\lambda-D^2)^{-1} \, \A^2_{[>0]} \cdots \A^2_{[>0]}\, 
 (\lambda-D^2)^{-1}\,  \A^2_{[>0]}\,(\lambda-D^2)^{-1}\right)_{[K]}\,  d\lambda\\
&=&  \sum_{j=0}^K   \sum_{\vert k\vert\geq 0}c(k_1, \cdots, k_j) \left(
\left(\A^2_{[>0]}\right)^{(k_1)}\left(\A^2_{[>0]}\right)^{(k_2)}\cdots
\left(\A^2_{[>0]}\right)^{(k_j)}\,\frac{1}{2i\pi} \int_{\Gamma} \,
\lambda^{\frac{2i-1}{2}} \,(\lambda-Q)^{-\vert k\vert-j-1}\, d\lambda
\right)_{[K]}\\
&=&  \sum_{j=0}^K   \sum_{\vert k\vert\geq
  0}\frac{\left(\frac{2i-1}{2}\right)\cdots\left(\frac{2i-1}{2}-\vert k\vert -j\right) }{(\vert k\vert +j)! }c(k_1, \cdots, k_j)\left(
\left(\A^2_{[>0]}\right)^{(k_1)}\left(\A^2_{[>0]}\right)^{(k_2)}\cdots
\left(\A^2_{[>0]} \right)^{(k_j)}\, \left(D^2\right)^{ 2i-1 -\vert k\vert-j}
\right)_{[K]}\\
\end{eqnarray*}
where the last equality follows by integration by parts. 
Applying the Wodzicki residue yields the result of the proposition. \endsquare

\section{ Getzler's rescaling}
 Let  ${\cal E}\to B$ be a \pdo vector bundle and let $\{D_b, b\in B\}$ be a smooth family of elliptic differential operators parametrised by $B$ acting on the fibres of 
${\cal E}$. \\
  Getzler's rescaling transforms  a homogeneous form $\alpha_{[i]}$ of degree $i$ to the expression
$$\delta_\e \cdot \alpha_{[i]}\cdot \delta_\e^{-1}=  \frac{\alpha_{[i]}}{ \sqrt \e^i},$$ so that a superconnection
$\A=\A_{[0]}+\A_{[1]}+\A_{[2]} $ transforms to 
$$\tilde \A_\e= \delta_\e \cdot \A\cdot \delta_\e^{-1}=\A_{[0]} +
\frac{\A_{[1]}}{\sqrt \e}+\frac{\A_{[2]}}{\e}. $$
Here we allow higher forms in the superconnection, keeping in mind later applications involving the Bismut superconnection for families of Dirac operators.  Following
 the usual conventions,  for a given superconnection $\A= \A_{[0]} +
\A_{[1]}+\A_{[2]}$  we set: 
$$ \quad \A_\e:= \sqrt \e\,\tilde \A_\e=  \sqrt \e\, \delta_\e
\cdot \A\cdot \delta_\e^{-1}=\sqrt \e \, \A_{[0]} +
\ \A_{[1]}+\frac{\A_{[2]}}{\sqrt \e}.$$
\begin{rk} Different notations are used in the literature, namely some authors
  set $t:=\sqrt \e$ which leads to (see e.g.\cite{L})  
  \begin{equation}\label{eq:LottsAt}
\bar   \A_t:=
  \A_{t^2}= t \, \A_{[0]} +
\ \A_{[1]}+\frac{\A_{[2]}}{t},
\end{equation}
a notation which we shall also use in this paper. 
\end{rk}
 The following result shows how the $j$-th Chern (resp. residue-) Chern form picks up the
 $2j$ (resp. $2j-1$-) form degree part of ${\rm str}^{ \A^2} \left(  \A^{2j} \right)$ . 
\begin{prop}\label{prop:resChernAt} Let $\A$ be a superconnection associated
  with a  family of  elliptic differential operators parametrised by $B$. Then, with the
  notations of (\ref{eq:LottsAt}) 
$${\rm str}^{ \A^2} \left(  \A^{2j} \right)_{[2j]}= {\rm fp}_{t=0} {\rm str} \left(
\bar \A_t^{2j} \, e^{-\bar \A_t^2}\right)_{[2j]}$$
and
$${\rm sres}\left(  \vert \A\vert ^{2j-1} \right)_{[2j-1]}=\frac{1}{2\sqrt \pi}  {\rm fp}_{t=0} {\rm str} \left(
\bar \A_t^{2j} \, e^{-\bar \A_t^2}\right)_{[2j-1]}.$$
\end{prop}
{\bf Proof:} Recall that since $D$ is a differential operator, so is $\A^2$ a
differential operator valued form and  (see Proposition\ref{prop:trQMellin})
$${\rm str}^{\A^2} \left(   \A^{2j}  \right)= {\rm fp}_{\e\to 0} {\rm str}
\left(   \A^{2j} e^{- \e \A^2} \right).$$
On the other hand, for any $t>0$ we have:
\begin{eqnarray*}
{\rm str} \left( \bar  \A_t^{2j}\,  e^{- \bar \A_t^2} \right)_{[2j]}&=& 
{\rm str} \left(  \A_{t^2}^{2j}\,  e^{-  \A_{t^2}^2} \right)_{[2j]}\\
&=& 
{\rm str} \left( t^{2j}\, \tilde \A_{t^2}^{2j}\,  e^{- t^2\, \tilde \A_{t^2}^2}
  \right)_{[2j]}\\
&=& 
{\rm str} \left( t^{2j} \, \delta_{t^2} \,  \A^{2j}\,\delta_{t^2}^{-1}\,  e^{-  t^2\,\delta_{t^2} \, \A^2\delta_{t^2}^{-1}}
  \right)_{[2j]}\\
&=& 
{\rm str} \left( t^{2j} \, \delta_{t^2} \, \left( \A^{2j}\,  e^{-
      t^2\,\A^2}\right) \delta_{t^2}^{-1} \right)_{[2j]}\\
&=& 
{\rm str} \left(  \A^{2j}\,  e^{-      t^2\,\A^2} \right)_{[2j]}.\\
 \end{eqnarray*}
Since the r.h. side is of the type ${\rm str} \left(  \A^{2j}\,  e^{-
    \e\,\A^2} \right)$ with $\A^2$  an elliptic differential operator
valued form  and $\A^{2j}\in \Omega(B, \Cl({\cal E}))$, it  has a
known  asymptotic expansion at $0$ and taking finite parts when $t\to 0$ yields the first part of the proposition.
\\ \\
Similarly, since $D$ is a differential operator we have $${\rm sres}\left(\vert \A\vert^{2j-1}\right)=
\frac{1}{2\sqrt \pi} {\rm fp}_{t=0} \left( t\, {\rm str} \left(\A^{2i} \,
    e^{-t^2\A^2}\right)\right)$$ and for any $t>0$
\begin{eqnarray*}
{\rm str} \left( \bar  \A_t^{2j} \, e^{- \bar \A_t^2} \right)_{[2j-1]}&=& 
{\rm str} \left(  \A_{t^2}^{2j} \, e^{-  \A_{t^2}^2} \right)_{[2j-1]}\\
&=& 
{\rm str} \left( t^{2j} \, \tilde \A_{t^2}^{2j} \,e^{- t^2\, \tilde \A_{t^2}^2}
  \right)_{[2j-1]}\\
&=& 
{\rm str} \left( t^{2j} \, \delta_{t^2}  \A^{2j}\, \delta_{t^2}^{-1}\,  e^{-
      t^2\,\delta_{t^2} \, \A^2\delta_{t^2}^{-1}}
  \right)_{[2j-1]}\\
&=& 
{\rm str} \left( t^{2j} \, \delta_{t^2}\,  \left( \A^{2j}\,  e^{-
      t^2\,\A^2}\right)
  \delta_{t^2}^{-1} \right)_{[2j-1]}\\
&=& t\,
{\rm str} \left(   \A^{2j}\,  e^{-
      t^2\,\A^2}\right)_{[2j-1]}.\\
\end{eqnarray*}
As before, since the r.h.s. is of the type ${\rm str} \left(  \A^{2j}\,  e^{-
    \e\,\A^2} \right)$ with  $\A^2$ an elliptic differential operator
valued form  and $\A^{2j}\in \Omega(B, \Cl({\cal E}))$, it  has a
known  asymptotic expansion at $0$ and taking finite parts when $t\to 0$ yields the first part of the lemma.
Taking finite parts when $t\to 0$ therefore yields the second part of the lemma. 
\endsquare\\ \\

\section{Superconnections associated with Dirac operators}

 We now specialise to Chern forms built from superconnections associated  with
 families of Dirac  operators. Let  $\pi:\M\to B$ be a smooth fibration of
 closed spin manifolds with fibre $M$ and $\E\to B$ a Clifford bundle  with an
 associated  family of Dirac operators parametrised by $B$. 
 The vector bundle ${\cal E}:= \pi_* \E$ is an (infinite rank) \pdo vector
 bundle with fibres modelled on $\Ci(\M/B, \E_{\M/B})$.
According to whether the manifolds are even or odd dimensional, ${\cal E}$ will be $\Z_2$-graded or not.\\ The  vertical Riemannian metric $g^{TM}$ and  hermitian 
metric $h^\E$  on $\E$ induce an $L^2$-inner product on ${\cal E}$. From a connection $\nabla^\E$ on $\E$ compatible with $h^\E$, one can build   a  connection  
$\tilde \nabla^{\cal E}_X \sigma:=\nabla^\E_{\tilde X} \sigma(b)$ on ${\cal E}=\pi_* \E$
where $\tilde X$ is the horizontal lift of $X\in T_bB$ and from there a unitary connection $\nabla^{\cal E}$ on ${\cal E}$.   \\
The corresponding Bismut superconnection associated with this fibration reads:
$$\A= D+ \nabla^{\cal E}+ c(T) \quad {\rm even }\quad {\rm case}, \quad
\A=\sigma \, D+ \nabla^{\cal E}+ \sigma\, c(T) \quad {\rm odd }\quad {\rm case},$$
where $T\in \Omega^2(\M, TM)$ is the curvature of the horizontal distribution on $\M$ and $c$ the Clifford multiplication. 
\\ Along the lines of the heat-kernel proof of the index theorem  we introduce the kernel $k_\e ( \A^{2})$ 
 of $e^{-\e \A^2}$ for some $\e>0$. Since $D$ is a family of Dirac
   operators, we have   (see e.g. chap. 10 in \cite{BGV}
 in the even dimensional case  and \cite{BC} in the odd dimensional case) 
\begin{equation} \label{eq:heatkernelexp}
k_\e(\A^{2})(x, x)\sim_{\e \to 0} \frac{1}{(4\pi \e)^\frac{n}{2}} \sum_{j=0}^\infty \e^j k_j( \A^{2})(x,
x)
\end{equation}\begin{prop}\label{prop:AChernforms} Let $\A$  be a superconnection adapted to a
  smooth family of Dirac  operators  parametrised by $B$.
\begin{enumerate}
\item 
  The $j$-th Chern form associated with $\A$ is given by an integration along
  fiber of $\M:$ 
 $${\rm str}^{\A^2} \left( \A^{2j}\right)=  \frac{(-1)^j j!}{(4\pi )^{\frac{n}{2}}} \int_{\M/B}{\rm str}( k_{j+\frac{n}{2}}(
      \A^{2})).$$
\item If  the kernel of $D$ has constant dimension,    the $j$-th residue Chern form associated with $\A$ reads: 
$$ \sqrt
\pi\, 
{\rm sres}  \,  \left(\vert \A\vert^{2j-1} \right)
= \frac{(-1)^j (2j-1)!!}{(4\pi )^{\frac{n}{2}} 2^{j-1}} \int_{\M/B}{\rm str}( k_{j+\frac{n-1}{2}}(
      \A^{2})). $$
Here  $\vert \A\vert^{2j-1}$ is defined as before by the
contour integral  $$\vert \A\vert^{2j-1} = \frac{i}{2\pi}\int_\Gamma
\l^{\frac{2j-1}{2}}\, (\A^2+\pi_\A - \lambda)^{-1} \; d\lambda.$$  
\end{enumerate}
\end{prop}
\begin{rk} It follows from this last formula   that adding a
    smoothing \pdo valued zero form to
  $\A^2$ does not affect the residue form ${\rm sres}  \,  \left(\vert
    \A\vert^{2j-1}\right)_{[2j-1]}$ so that one expects formulae for the
  residue Chern form to be independent of $\pi_\A$.
\end{rk} 
{\bf Proof:} The trace under the integral sign is just the matrix trace for endomorphims of a finite rank vector bundle $E$ whereas 
on the left-hand-side the trace is computed in the Hilbert space of
square-integrable sections of $E$ over $M.$\\ \\ By the above remark,   $\A$ being  a differential operator, it has vanishing Wodzicki residue and we have:
$${\rm str}^{\A^2} \left( \A^{2j}\right)=   {\rm fp}_{\e = 0} {\rm str}\left(
\A^{2j} \, e^{-\e \A^2}\right)=(-1)^j
{\rm fp}_{\e = 0}\left(\partial_\e^j {\rm str}\left(e^{-\e \A^2}\right)_{\vert_{\e= 0}}\right).$$
 By equation (\ref{eq:heatkernelexp}), this yields
\begin{eqnarray*}
{\rm str}^{ \A^2} \left( \A^{2j}\right)&= & (-1)^j
{\rm fp}_{{\e= 0}}\left(\partial_\e^j {\rm str}\left(e^{-\e\,   \A^2}\right)\right)\\
&=& (-1)^j\, {\rm fp}_{{\e= 0}}\left(
\partial_\e^j \int_{\M/B}{\rm str}\left(k_\e( \A^{2})\right)\right)\\
&=  &    \frac{(-1)^j}{(4\pi )^{\frac{n}{2}} }
\, {\rm fp}_{{\e= 0}}\left(\partial_\e^i \sum_{j=0}^\infty \e^{j-\frac{n}{2}} \int_{\M/B} {\rm str}( k_j(\A^{2}))\right)\\
&= &    \frac{(-1)^i}{(4\pi )^{\frac{n}{2}}} 
\, {\rm fp}_{{\e= 0}}\left( \sum_{l=0}^\infty (l -\frac{n}{2})\cdots (l -\frac{n}{2}-j+1)\e^{l -\frac{n}{2}-j}\int_{\M/B} {\rm str}( k_l(
      \A^{2})\right)\\
&=&  \frac{(-1)^j j!}{(4\pi )^{\frac{n}{2}}} \int_{\M/B}{\rm str}( k_{j+\frac{n}{2}}(
      \A^{2})).\\
\end{eqnarray*}
This proves the first part of the proposition. 
On the other hand,  again by (\ref{eq:heatkernelexp}) we have:
\begin{eqnarray*}
{\rm fp}_{\e = 0}{\rm str}\left(\sqrt \e \,\left(\A^{2j}e^{-\e \,\A^2}\right)\right)&=&
(-1)^j{\rm fp}_{{\e= 0}}\left(\sqrt \e\,\,  \partial_\e^j {\rm str}\left(e^{-\e\,   \A^2}\right)\right)\\
&=& (-1)^j, {\rm fp}_{{\e= 0}}\left(\sqrt \e\, 
\partial_\e^j \int_{\M/B}{\rm str}\left(k_\e( \A^{2})\right)\right)\\
&=  &    \frac{(-1)^j}{(4\pi )^{\frac{n}{2}}} 
\, {\rm fp}_{\e= 0}\left(\sqrt \e\,\partial_\e^j \sum_{l=0}^\infty \e^{l-\frac{n}{2}} \int_{\M/B} {\rm str}( k_l(
      \A^{2}))\right)\\
&= &    \frac{(-1)^j}{(4\pi )^{\frac{n}{2}}} 
\, {\rm fp}_{{\e= 0}}\left( \sum_{l=0}^\infty (l -\frac{n}{2})\cdots (l -\frac{n}{2}-j+1)\e^{l -\frac{n-1}{2}-j}\int_{\M/B} {\rm str}( k_l(
      \A^{2}))\right)\\
&=&  \frac{(-1)^j (j-\frac{1}{2}) (j-\frac{3}{2})\cdots \frac{1}{2}}{(4\pi )^{\frac{n}{2}}} \int_{\M/B}{\rm str}( k_{j+\frac{n-1}{2}}(
      \A^{2}))\\
&=& \frac{(-1)^j (2j-1)!!}{(4\pi )^{\frac{n}{2}} 2^{j-1}}\, \int_{\M/B}{\rm str}( k_{j+\frac{n-1}{2}}(
      \A^{2})).
\end{eqnarray*}
  On the other hand,  Lemma \ref{lem:Mellin} applied  to $g(\e)= {\rm str} \left(\A^{2j} e^{-\e \A^2}\right)$  then  yields:
\begin{eqnarray*}
 {\rm fp}_{\e =0}\left( \sqrt \e \, {\rm str} \left(\A^{2j} e^{-\e \A^2}\right)\right)&=&  \Gamma(\frac{1}{2})\,{\rm res}_{z=0}  \, {\rm str} \left(\A^{2j} 
(\A^2)^{-z-\frac{1}{2}}\right) \\
&=& 2\, \sqrt \pi \, {\rm sres}  \,  \left(\A^{2j} (\A^2)^{-\frac{1}{2}}\right) \\
\end{eqnarray*}
 so that, 
$$ \sqrt \pi\, {\rm sres}  \,  \left(\A^{2j}
  \left(\A^2 \right)^{-\frac{1}{2}}\right)= \frac{(-1)^j (2j-1)!!}{(4\pi )^{\frac{n}{2}} 2^{j}}\, \int_{\M/B}{\rm str}( k_{j+\frac{n-1}{2}}( \A^{2}))$$

which  proves the second part of the  the proposition. 
\endsquare\\ \\
The following result relates the $j$-th  (resp. residue) Chern   form with the
$2j$-th (resp. $2j-1$-th) form degree part of the Chern character $\lim_{t\to 0}{\rm ch}(\A_t)$.
\begin{thm}\label{thm:localChernforms} In the $\Z_2$-graded  case the $j$-th Chern form associated with
  a superconnection $\A$  reads: 
 \begin{eqnarray}\label{eq:strChernj}
{\rm str}^{\A^2} \left( \A^{2j}\right)_{[2j]}&=& \frac{(-1)^j j!}{(2i\pi)^{\frac{n}{2}}} \, \left( \int_{\M/B} \hat A(\M/B) \wedge
  {\rm ch}(\E_{\M/B}) \right)_{[2j]}\nonumber\\
&=& (-1)^j j! \,\left( \lim_{t\to 0}{\rm ch}(\A_t) \right)_{[2j]}.
\end{eqnarray}
In the ungraded  case the $j$-th residue Chern form associated with the
superconnection with kernel of $D$ of constant dimension reads: 
\begin{eqnarray}\label{eq:resChernj}
 {\rm sres}  \,  \left(\vert \A\vert^{2j-1}\right)_{[2j-1]}
&=&  \frac{(-1)^j (2j-1)!!}{(2i\pi)^{\frac{n+1}{2}}\, 2^{j-1}} \, \left(\int_{\M/B} \hat A(\M/B) \wedge
  {\rm ch}(\E_{\M/B}))\right)_{[2j-1]}\nonumber\\
&=&    \frac{(-1)^j (2j-1)!!}{ 2^{j-1}\sqrt \pi}  \,  \left( \lim_{t\to 0}{\rm ch}(\A_t) \right)_{[2j-1]}. 
\end{eqnarray}
\end{thm}
{\bf Proof:} 
As in   \cite{BGV} par. 10.4, using
the asymptotic expansion  of the  kernel $k_t(x, x)$ of the heat-operator
$e^{- t \A^2}$:
$$k_t(x,x)\sim_{t\to 0}\frac{1}{(4\pi t)^{\frac{n}{2}}}\sum_{j=0}^\infty t^j\,
  k_j(x, x)$$
we have:
\begin{eqnarray*}
{\rm ch}(\A_t)&=& \delta_t \left( {\rm str} (e^{- t \A^2})\right)\\
&\sim_{t \to 0}& (4 \pi t)^{-\frac{n}{2}} \sum_j t^j\int_{\M/B} \delta_t \left({\rm str}( k_j (\A^2))\right)\\
&\sim_{t \to 0} & (4 \pi )^{-\frac{n}{2}} \sum_{j, p} t^{j-(n+p)/2}\left( \int_{\M/B}{\rm str} \left(   k_j (\A^2)\right)\right)_{[p]},\\
\end{eqnarray*} 
so that 
\begin{equation}\label{eq:rescaledkernelasymp}
{\rm fp}_{t=0}{\rm ch}(\A_t)_{[p]}= (4 \pi )^{-\frac{n}{2}}  \left(
  \int_{\M/B}{\rm str} \left(   k_{\frac{p+n}{2}} (\A^2)\right)\right)_{[p]}.
\end{equation}
\begin{itemize}
\item {\bf  $\Z_2$-graded case. }
The family index  theorem  \cite{B} (see also  Theorem 10.23 in \cite{BGV})
yields the existence of the  limit as $t\to 0$ and 
$$\lim_{t\to  0} {\rm ch}(\A_t)= (2 i\pi)^{-\frac{n}{2}} \, \int_{\M/B} \hat A (\M/B)\wedge
  {\rm ch}(\E_{\M/B}).$$
Combining these two facts leads to:
$$\left( \int_{\M/B}{\rm str} \left(   k_{\frac{n+2j}{2}} (\A^2)\right)\right)_{[2j]}=\frac{(4 \pi )^{\frac{n}{2}}}{(2i\pi)^{\frac{n}{2}}} \left( \int_{\M/B} \hat A(\M/B) \wedge
  {\rm ch}(\E_{\M/B}) \right)_{[2j]}.$$
Inserting this in  Proposition \ref{prop:AChernforms} yields $$ {\rm str}^{\A^2} \left( \A^{2j}\right)_{[2j]}=  \frac{(-1)^j j!}{(2i\pi)^{\frac{n}{2}}} \, \left(\int_{\M/B} \hat A(\M/B) \wedge
  {\rm ch}(\E_{\M/B})\right)_{[2j]}$$ and hence (\ref{eq:strChernj}). 
\item {\bf  Ungraded  case. }
 The family index  theorem   yields the existence of the
  limit and             \cite{BC}.
$$\lim_{t\to 0} {\rm ch}(\A_t)=\frac{\sqrt \pi}{  (2\pi i)^{\frac{n+1}{2}}} \, \int_{\M/B} \hat A (\M/B)\wedge
  {\rm ch}(\E_{\M/B}).$$
 Combined with (\ref{eq:rescaledkernelasymp})
this   yields:
$$\left( \int_{\M/B}{\rm str} \left(   k_{\frac{n+2j-1}{2}}
    (\A^2)\right)\right)_{[2j-1]}=\sqrt \pi\,   \frac{(4 \pi )^{\frac{n}{2}}}{(2i\pi)^{\frac{n+1}{2}}}\, \left(\int_{\M/B} \hat A(\M/B) \wedge
  {\rm ch}(\E_{\M/B})\right)_{[2j-1]}.$$
Inserting this in Proposition \ref{prop:AChernforms} gives:
$$
 {\rm sres}  \,  \left(\vert \A\vert^{2j-1}\right)_{[2j-1]}
=  \frac{(-1)^j (2j-1)!!}{(2i\pi)^{\frac{n+1}{2}}\, 2^{j-1}}\,  \left(\int_{\M/B} \hat A(\M/B) \wedge
  {\rm ch}(\E_{\M/B}))\right)_{[2j-1]} $$
and hence (\ref{eq:resChernj}).
 \end{itemize}
\endsquare \\
\\
\begin{cor}
Whenever the fibration $\M\to B$ is
  trivial, then
\begin{equation}\label{eq:splittingA} 
 {\rm tr}^Q(\Omega^j)- \frac{(-1)^j j!}{(2i\pi)^{\frac{n}{2}}} \,\left(\int_{\M/B} \hat A (\M/B)\wedge
  {\rm ch}(\E)\right)_{[2j]}=\bar f_j(\bar \nabla).
\end{equation}
is local in the sense of the above definition.
\end{cor} 
\begin{ex} In the $\Z_2$-graded  case,  $ \frac{1}{(2i\pi )^{\frac{n}{2}}} \, \left(\int_{\M/B} \hat A (\M/B)\wedge
  {\rm ch}(\E_{\M/B})\right)_{[2]}$
corresponds to the curvature on the determinant line bundle associated with
a family of Dirac operators \cite{BF}. The formula corresponding to $j=1$ in the above theorem 
$${\rm str}^{\A^2} \left( \A^{2}\right)_{[2]}= - \frac{1}{(2i\, \pi )^{\frac{n}{2}}}\,  \left(\int_{\M/B} \hat A (\M/B)\wedge
  {\rm ch}(\E_{\M/B})\right)_{[2]}$$
expresses  the curvature on the determinant bundle  as $- \frac{1}{(2i
  )^{\frac{n}{2}}}$ times  the degree $2$ part of the first Chern form
associated with the superconnection $\A$, thereby generalising  the relation
that holds in finite dimensions (corresponding to the case $n=0$ of a zero
dimensional fibre $M$) relating the first Chern form on a finite rank
supervector bundle with minus the curvature on its determinant bundle (see \cite{PR}
for a discussion concerning this relation).
\end{ex}
\begin{ex}In the ungraded case, $  \left(\int_{\M/B} \hat A \wedge
  {\rm ch}(\E_{\M/B}))\right)_{[2j-1]}$  corresponds to the curvature of a
gerbe with connection associated with the family of Dirac operators \cite{CM},\cite{EM},\cite{L}.
The formula obtained in the above theorem for $j=2$
\begin{equation}\label{eq:gerbecurv}
 {\rm sres}  \,  \left(\A^2\, \vert \A\vert \right)_{[3]}
= \frac{3 }{2\, (2i\pi)^{\frac{n+1}{2}}} \, \left(\int_{\M/B} \hat A \wedge
  {\rm ch}(\E_{\M/B}))\right)_{[3]} 
\end{equation}
where we have set $\vert \A \vert=
   (\A^2)^{\frac{1}{2}}$
 relates this curvature with the degree $3$ part of the residue Chern form ${\rm sres}  \,  \left(\A^2\, \vert \A\vert \right)_{[3]}$.
\end{ex}
\section{Transgressed residue  Chern forms}
 Let  as before  $\pi:\M\to B$ be a smooth fibration of
 closed odd dimensional  spin manifolds with fibre $M$ and $\E\to B$ a Clifford bundle. 
 The vector bundle ${\cal E}:= \pi_* \E$ is an (infinite rank) \pdo vector
 bundle with fibres modelled on $\Ci(\M/B, \E_{\M/B})$.
Let  $\{D_b, b\in B\}$  be a smooth family of   Dirac  operators
associated with this fibration.\\
We  need to work with invertible  operators and introduce for
this purpose a  covering of $B$ by  open sets   $\{U_\lambda\}_{\lambda\in \R}$ with the
property that for any $b\in U_\lambda$  the (discrete) spectrum of $D_b$ does
not contain  $\lambda$. Then  $D(\lambda)= D- \lambda I$ is a family of Dirac
type operators which is everywhere invertible on
$U_\lambda$ and
$$ \A_{\lambda}:=  \sigma \, D(\lambda)+ \nabla^{{\cal E}}+ 
\sigma \, \frac{c(T)}{4}$$ is a  superconnection
associated with $D(\lambda)$.\\

 With the notations of (\ref{eq:LottsAt}) we set
$$\bar \A_{t}:= t \, \sigma \, D+ \nabla^{{\cal E}}+ \sigma\frac{ c(T)}{4t}
,$$ 
where $\sigma$ is the grading, which defines a smooth family of  superconnections
adapted to $D$.\\ 
Let for $t>0$ $$ \bar\A_{\lambda, t}:=  \sigma\, t\, D(\lambda)+ \nabla^{{\cal E}}+ 
\sigma \, \frac{c(T)}{4t},$$ 
 which defines a smooth family of  superconnections
adapted to $D(\lambda)$.
 The following technical result will be useful for
  what follows. 
\begin{lem}\label{lem:asympbehaviour} Let $\gamma$ be a differential
    operator valued form on $B$.
\begin{enumerate}
\item  The function 
$$t\mapsto {\rm str} \left( \gamma \, e^{-\bar \A_{\lambda,t}^2}\right)$$
decreases faster than any power of $t$  as $t\to \infty$. As     $t\to 0$ it
behaves as a  finite linear
combination of expressions 
 $\sum_{j=0}^\infty \alpha_j
t^{j-\delta}$
for some integer  $\delta$  and complex numbers $ \alpha_j, \beta_k$.
\item If  $\gamma$ is  an {\bf even} form, then 
\begin{enumerate}
\item $ {\rm str} \left( \sigma\, \gamma \, e^{-\bar
      \A_{\lambda,t}^2}\right)$ only involves odd powers of $t$ and $t^{-1}$. In
  particular, $${\rm fp}_{t=0} {\rm str} \left( \sigma\, \gamma \, e^{-\bar
      \A_{\lambda,t}^2}\right)=0.$$
\item For any non negative integer $j$
  $$\left({\rm str} \left( \sigma\, \gamma \, e^{-\bar
      \A_{\lambda,t}^2}\right)\right)_{[2j]}=0\quad \forall\quad  t>0.$$
\end{enumerate}
\end{enumerate} 
 \end{lem}
\begin{rk} The  lemma  easily extends to \pdo valued forms if we allow for logarithmic
  divergences in $t$ in which case $\delta$ is a real number. \end{rk}
 {\bf Proof:}\begin{enumerate}
\item
 Let us first introduce   notations similar to notations
  of \cite{H}.  Let $T\in \Cl(M, E)$ and $T_k, k\in \N$
be operators in $\Cl(M, E)$  with decreasing order in $k$. Then 
\begin{eqnarray}\label{eq:sim}
&{}& T\sim \sum_{k\geq 0} T_k\\
&\Longleftrightarrow&\exists C\in \Cl(M, E) \quad {\rm invertible}, {\rm
  s.t.}\quad  \forall N \in \N, \exists K(N)\nonumber\\
&{}&
\quad \left(T- \sum_{k= 0}^{K(N)} T_k\right)\, C\in \Cl^{-N}(M, E). \nonumber
\end{eqnarray} 
In the sequel, the operator $e^{-t^2 D^2}$  plays the role of the invertible operator $C$.
\\ We also need to   extend to \pdo valued forms,  notations
previously used for ordinary classical pseudo-differential operators.
 For    $\beta\in \Omega\left((B, \Cl\left({\cal E}\right)\right)$   and  any $j\in \N$ we set: $$\alpha^{(j)}(\beta):= {\rm ad}_{D^2}^j(\beta),
\quad {\rm where} \quad  {\rm ad}_{ D^2}(\beta)=[D^2, \beta],$$  so that $ \beta^{(0)}=\beta, \quad
 \beta^{(j+1)}= {\rm ad}_{D^2}(\beta^{(j)})=[D^2, \beta^{(j)}].$\\ \\
Since $\bar \A_{\lambda, t}= t\, \sigma \, D(\lambda) +\left(\bar \A_{\lambda, t}\right)_{[>0]}$ we 
have \begin{eqnarray*}
\bar  \A_{\lambda, t}^2&= &t^2\,D(\lambda)^2+\left(\bar  \A_{\lambda,
    t}^2\right)_{[>0]}\\
&=& t^2\,D(\lambda)^2+\sigma\,\left[ t\, [\nabla^{\cal E}, D(\lambda)]+ 
\frac{[\nabla^{\cal E}, c(T)]}{4t}\right]+ \frac{1}{4}
[D(\lambda), c(T)]+ \frac{c^2(T)}{16\, t^2} +\Omega^{\cal E}\\
&=& t^2\,D(\lambda)^2+\sigma\, \left(\bar \A_{\lambda, t}\right)_{[>0, od]}+\left(\bar \A_{\lambda, t}\right)_{[>0,ev]}
\end{eqnarray*}
where we have set  $\left(\bar \A_{\lambda, t}\right)_{[>0, od]}:= t\, [\nabla^{\cal E}, D(\lambda)]+ 
\frac{[\nabla^{\cal E}, c(T)]}{4t}$ which only involves odd powers of $t$ and $t^{-1}$.
\\  
  Duhamel's formula then yields (see e.g. \cite{H}):
\begin{eqnarray*}
   e^{ -\bar \A_{\lambda, t}^2}&=&
 (-1)^n\int_{\Delta_l}e^{- u_0\, t^2 D^2} \left(\bar \A_{\lambda,t}^2\right)_{[>0]} \cdots  e^{- u_{l-1} \,
    t^2 D(\lambda)^2}\, \left(\bar \A_{\lambda,t}^2\right)_{[>0]}\, e^{- u_l \, t^2 D(\lambda)^2}\, du_1\cdots \, du_l\\
&\sim &\sum_{\vert k\vert\geq 0}\frac{(-1)^{\vert k\vert} t^{2\vert k\vert}\, c(k)}{(\vert
  k\vert+n)!} \left(\left(\bar \A_{\lambda,t}^2\right)_{[>0]}\right)^{(k_1)} \cdots
 \left(\left(\bar \A_{\lambda,t}^2\right)_{[>0]}\right)^{(k_l)}e^{- t^2
   D(\lambda)^2}.\\
\end{eqnarray*} 
Here $\Delta_l:=\{(u_0, \cdots, u_l), u_i\geq 0, \quad \sum_{i=0}^l u_i=1\}$
is the unit simplex and with the coefficient $c(k)$ as previously defined.
\\ 
Since $\left(\bar \A_{\lambda,t}^2\right)_{[>0]}$ has positive degree,  only a finite number of terms of the sum will contribute for a fixed form
degree.\\
Now, for a differential operator  $C$ of order $c$, the map  $t\mapsto {\rm str} \left(C\,
  e^{-t^2 D^2}\right)$ decreases faster than any power of $t$  at
infinity and 
behaves asymptotically as follows  as $t\to 0$: 
\begin{equation}\label{eq:asymptotics}
{\rm str} \left(C\, e^{-t^2 D^2}\right)\sim_{t\to 0}\sum_{j=0}^\infty \alpha_j
t^{j-c-n} \, 
\end{equation}
 where $c$ is the order of $C$ and $n$ the dimension of the manifold.\\
Since  expressions of the type  $\left(\bar \A_{\lambda,t}^2\right)_{[>0]}$ are linear
combinations of differential operator 
valued forms with coefficients given by  powers of $t$, for any \pdo
valued form $\gamma$ on $B$
$${\rm str}\left( \gamma\, \left(\left(\bar \A_{t}^2\right)_{[>0]}\right)^{(k_1)} \cdots
 \left(\left(\bar \A_{t}^2\right)_{[>0]}\right)^{(k_n)}e^{- t^2 D^2} \right)
$$ has the expected asymptotic  behaviour at $0$ and at $\infty.$ It follows
that so does ${\rm str}\left( \gamma\, e^{ -\bar \A_{t}^2}\right)$ have a
similar asymptotic behaviour. This ends the proof of the first part of the
lemma.
\item The second part of the lemma requires a closer look at the expressions
  involved. Since
$$\left(\bar \A_{\lambda, t}\right)_{[>0]}=  \sigma\, \left(\bar \A_{\lambda, t}\right)_{[>0, od]}+\left(\bar \A_{\lambda, t}\right)_{[>0,ev]}$$
and since  str vanishes on  terms of the type $\sigma
\beta$,  in the expression 
\begin{eqnarray*}
&{}&{\rm str} \left( \sigma \, \gamma\,  e^{ -\bar \A_{\lambda, t}^2}\right)\\
&=& \sum_{\vert k\vert\geq 0}\frac{(-1)^{\vert k\vert} t^{2\vert k\vert}\, c(k)}{(\vert
  k\vert+n)!}{\rm str} \left(  \sigma \, \gamma\,\left(\left(\bar \A_{\lambda,t}^2\right)_{[>0]}\right)^{(k_1)} \cdots
 \left(\left(\bar \A_{\lambda,t}^2\right)_{[>0]}\right)^{(k_l)}\, e^{- t^2
   D(\lambda)^2}\right)
\end{eqnarray*}
(which we recall only contains a finite number of terms for fixed form
degree)   only those terms will remain that involve an odd number of expressions of
the type $\left(\bar \A_{\lambda, t}\right)_{[>0, od]}$ and hence odd
  powers of $t$ and $t^{-1}$.  Since $\gamma$ is
assumed to be of even degree, it follows that  the
total expression ${\rm str} \left( \sigma \, \gamma\,  e^{ -\bar \A_{\lambda,
      t}^2}\right)$ is an odd degree form which only involves odd powers of $t$
and $t^{-1}$ so that $${\rm fp}_{t=0}\left({\rm str} \left( \sigma \, \gamma\,  e^{ -\bar \A_{\lambda,
      t}^2}\right)\right)=0 ;\quad \left({\rm str} \left( \sigma \, \gamma\,  e^{ -\bar \A_{\lambda,
      t}^2}\right)\right)_{[2j]}=0\quad \forall t>0.$$
 This proves the second part of the lemma.
\end{enumerate}
 \endsquare \\
\\

 The following theorem provides a  transgression formula for the residue  Chern forms  $ {\rm sres}\left( \vert
  \A\vert^{2j-1}\right)_{[2j-1 ]}$. 
\begin{thm}\label{thm:transgression} On every open subset $U_\lambda$, the
  following transgression formula holds:
$$
 {\rm sres} \left( \vert
 \A_\lambda\vert^{2j-1}\right)_{[2j-1 ]}=a_j \cdot  
  d\, \left(\tilde\eta_\lambda\right)_{[2j-2]},$$
where $d$ is the exterior differential, $a_j:= 
\frac{(-1)^j (2j-1)!!}{2^{j-1} } $ and $$\tilde \eta_\lambda = {\rm
   fp}_{t=0} \int_t^\infty {\rm str} \left[  
 \frac{d}{ds}\bar \A_{\lambda,s}\, e^{-\bar
 \A_{\lambda,s}^2}\right]\,ds$$
is the $\eta$ invariant associated with $D(\lambda)$ (see \cite{BC}, \cite{L}). 
\end{thm}
{\bf Proof:} We first show that  ${\rm sres}\left( \vert
  \A\vert^{2j-1}\right)_{[2j-1 ]}={\rm sres}\left( \vert
  \A_\lambda\vert^{2j-1}\right)_{[2j-1 ]}$ and then show a transgression
formula for ${\rm sres}\left( \vert
  \A_\lambda\vert^{2j-1}\right)_{[2j-1 ]}$. 
\begin{enumerate}
\item Let us consider a smooth family $D(\lambda)(\e):= D-\e \,\lambda
\, I$ of first order elliptic differential operators interpolating $D$ and
$D(\lambda)$ between $0$ and $1$ and the corresponding superconnections $$\bar
\A_{\lambda,t}(\e):= 
\A_{\lambda, t} -\e\,\sigma\, \lambda \, I.$$
Differentiating w.r. to $\e$ we have:
\begin{eqnarray*}
\frac{d}{d\e} {\rm str} \left(e^{-\bar \A^2_{\lambda, t} (\e)}\right)&=&
- {\rm str} \left([\bar \A_{\lambda,t}(\e),\frac{d}{d\e}\bar \A_{\lambda,t}(\e)]\,
  \,e^{-\bar \A^2_{\lambda, t} (\e)}\right)\\
&=&
-\, {\rm str} \left([\bar \A_{\lambda,t},\frac{d}{d\e}\bar \A_{\lambda,t}(\e)\,  \,e^{-\bar \A_{\lambda, 
       t}^2(\e)}]\right)\\ 
&=&- d\, {\rm str} \left(\frac{d}{d\e}\bar \A_{\lambda,t}(\e)  \,e^{-\bar \A_{\lambda, 
       t}^2(\e)}\right)\\ 
&=&\lambda\, d\,  {\rm str} \left( \sigma \,e^{-\bar \A_{\lambda, 
       t}^2(\e)}\right).\\ 
 \end{eqnarray*}
By part 2 of Lemma \ref{lem:asympbehaviour} applied to $\gamma= I$, we find
that  for any positive  integer $j$
$$\frac{d}{d\e} {\rm str} \left(e^{-\bar \A^2_{\lambda, t}
    (\e)}\right)_{[2j-1]}=\lambda\, d\, \left[ {\rm str} \left( \sigma \,e^{-\bar \A_{\lambda, 
       t}^2(\e)}\right)\right]_{[2j-2]}=0$$
as a consequence of which ${\rm str} \left(e^{-\bar \A^2_{\lambda, t}
    (\e)}\right)_{[2j-1]}$ is actually  independent of $\e$ and 
 $${\rm str} \left(e^{-\bar \A^2_{\lambda, t}
   }\right)_{[2j-1]}=  {\rm str} \left(e^{-\bar \A^2_{ t}
    }\right)_{[2j-1]}.$$
By formula (\ref{eq:resChernj}) in Theorem \ref{thm:localChernforms}, the limit
on either side therefore exists as $t\to 0$ and  
\begin{eqnarray}\label{eq:Atlambda}
\lim_{t\to 0}{\rm str} \left(e^{-\bar \A^2_{\lambda, t} }\right)_{[2j-1]}
&=& \lim_{t\to 0} {\rm ch}\left(\A_t\right)_{[2j-1]}\nonumber\\
&=&
   \frac{(-1)^j2^{j-1}}{ (2j-1)!! }   \,  {\rm sres}\left( \vert
  \A\vert^{2j-1}\right)_{[2j-1 ]}.\end{eqnarray}
\item We now derive a transgression formula for ${\rm str} \left(e^{-\bar \A^2_{\lambda, t}
    }\right)_{[2j-1]}$, from which will then follow a transgression formula
  for ${\rm sres}\left( \vert
  \A\vert^{2j-1}\right)_{[2j-1 ]}$ as a consequence of (\ref{eq:Atlambda}).

\begin{eqnarray}\label{eq:derivative}
\frac{d}{dt} {\rm str} \left(
 e^{-\bar \A_{\lambda,t}^2}\right)&=& \,  - {\rm str} \left(
 [\bar \A_{\lambda,t}, \frac{d}{dt}\bar \A_{\lambda,t}]\, e^{-\bar \A_{\lambda,t}^2}\right)\nonumber\\
&=& - {\rm str} \left([ \bar \A_{\lambda,t},
 \,  \frac{d}{dt}\bar \A_t\, e^{-\bar \A_{\lambda,t}^2}]\right)\nonumber\\
&=&  - d\, {\rm str} \left(
 \, \frac{d}{dt}\bar \A_{\lambda,t}\, e^{-\bar \A_{\lambda,t}^2}\right).
\end{eqnarray}
 The first part of  Lemma \ref{lem:asympbehaviour}
 provides  a control   as $t\to 0$ and as 
 $t\to \infty$ on the asymptotic behaviour of the last expression
$ {\rm str} \left[
\, \frac{d}{ds}\bar \A_{\lambda,s}\, e^{-\bar \A_{\lambda,s}^2}\right] $ arising in (\ref{eq:derivative}).
 Its  primitive in $t$ 
$$\tilde \eta_\lambda(t):= \int_t^\infty {\rm str} \left[    \frac{d}{ds}\bar \A_{\lambda,s}\, e^{-\bar
 \A_{\lambda,s}^2}\right]\, ds$$ which exists as a
  consequence of the invertibility of $D(\lambda)$,  has a similar asymptotic
behaviour as $t\to 0$. Integrating (\ref{eq:derivative}) from $t$ to $\infty$ and taking the finite part as $t\to 0$  we find that  the $\eta$
invariant (we borrow notations from \cite{L}, see his formula (3.19))
 $$\tilde \eta_\lambda:= {\rm fp}_{t= 0} \tilde \eta_\lambda(t) = {\rm
   fp}_{t=0} \int_t^\infty {\rm str} \left[  
 \frac{d}{ds}\bar \A_{\lambda,s}\, e^{-\bar
 \A_{\lambda,s}^2}\right]ds$$
trangresses  ${\rm fp}_{t=0}{\rm str} \left(e^{-\bar \A^2_{\lambda, t}}\right)$:
$${\rm fp}_{t=0}{\rm str} \left(e^{-\bar \A^2_{\lambda, t }}\right)         = d\, \tilde\eta_\lambda.$$
But by the first part of the theorem (see equation (\ref{eq:Atlambda})), this leads to 
$$
 {\rm sres} \left( \vert
 \A_\lambda\vert^{2j-1}\right)_{[2j-1 ]}= \frac{(-1)^j (2j-1)!!}{2^{j-1} }
  d\,\left( \tilde\eta_\lambda\right)_{[2j-2 ]}.$$
\end{enumerate}

\section{Relation to hamiltonian  gauge anomalies}
We  first review a finite dimensional situation which will serve as a
model for infinite dimensional genralisations. We consider  the finite-dimensional Grassmann manifold ${\rm Gr}(n,n)$
consisting of rank $n$ projections in
$\Bbb C^{2n}$, which we   parametrise by grading
operators $F= 2P -1,$ where $P$ is a finite rank projection.
\begin{lem} The even forms 
\begin{equation}\label{eq:finitedimomegaj} \omega_{2j} = \tr\left( \, F(dF)^{2j}\right), \end{equation}
where $j=1,2,\dots$ are closed  forms on ${\rm Gr}(n,n)$.
 \end{lem}
{\bf Proof:} By the traciality of $\tr$ we have
\begin{eqnarray}
 d\, \omega_{2j} &=& d\, \tr\left( \, F(dF)^{2j}\right)\nonumber\\
&=& \tr\left( (dF)^{2j+1}\right)\nonumber\\
&=& \tr\left( F^2\, (dF)^{2j+1}\right)\nonumber\\
&{}&{\rm since} \quad F^2=1\nonumber\\
&=& -\tr\left( F\, (dF)^{2j+1}\, F\right)\nonumber\\
&{}&{\rm since} \quad F\, dF=-dF \, F\nonumber\\
&=& -\tr\left( (dF)^{2j+1}\, F^2\right)\nonumber\\
&{}&{\rm since} \quad \tr([A, B])=0\nonumber\\
&=& -\tr\left( (dF)^{2j+1}\right)\nonumber\\
&=& 0.
\end{eqnarray}
\endsquare\\ \\
In fact it turns out that the cohomology of  ${\rm Gr}(n, n)$ is 
generated by even (nonnormalized) forms of the type $\omega_{2j}, j=1, \cdots,
n$ \cite{MS}.
\\ \\
 Let us now consider the  infinite dimensional geometric setup described  in
 the previous section up to the fact that $\pi: \M=M\times B \to B$ is now a trivial
 fibration with typical fibre a closed (Riemannian) spin manifold $M$.
\\ On each open subset $U_\lambda:=\{b \in B,\lambda\notin {\rm
  spec}(D_b)\}\,  \subset B$ there is a well defined map
\footnote{Note here again,
the difference in conventions compared to \cite{L}. We follow \cite{CMM}, \cite{CM}, whereas in \cite{L} the operators $D_b$ are perturbed as $D_b\mapsto D_b +h(D_b)$ by 
a smoothing function $h$ to avoid the zero modes.} 
\begin{eqnarray*}
F: B&\to & \Cl_0(M, E)\\
b &\mapsto &
F_b:=(D_b-\lambda I)/|D_b -\lambda I|.
\end{eqnarray*} 
 Since $F_b^2= F_b$, $P_b:= \frac{I+ F_b}{2}$ is a projection, the
  range Gr$(M, E):= {\rm Im} F$ of $F$ coincides with   the Grassmannian consisting of classical pseudodifferential  projections $P$ with kernel and
cokernel of infinite rank, acting in the complex Hilbert space $H:= L²(M,
E)$. 
Here $L^2(M, E)$ denotes the space of  square-integrable sections of the vector bundle $E$ over the
compact manifold $M.$\\
  This map $b\mapsto F_b$  is generally not contractible and we want
to define cohomology classes on $B$ as in (\ref{eq:strChernj}) up to some
modifications required by the specific situation.  This problem usually arises 
in hamiltonian quantization in field theory, when $M$ is an odd dimensional manifold, the physical 
space. Although here we deal with even forms for odd order  operators, there is a relation to the previous discussion 
on odd forms for odd order operators which is explained in the end of this section.
The problem here 
is similar in spirit to the earlier discussion in as far as  we want to modify the naive cohomology classes,
imitating the finite-dimensional case, by local corrections arising from the infinite dimensionality of 
the problem. 
\\ 
Indeed, in this infinite dimensional setup  traces are  generally 
ill-defined, so that we cannot a priori extend the above computation to 
${\rm Gr}(M, E)$. \footnote{Unless  we restrict to the submanifold ${\rm Gr}_{res}(M, E) \subset {\rm Gr}(M, E)$ consisting of points $F$ such that 
$F-\epsilon$ is Hilbert-Schmidt for some fixed point $\epsilon \in {\rm Gr}(M,
E)$ \cite{PS}}. 
\\ However, we can define an analog of  (\ref{eq:finitedimomegaj}) at the cost
of replacing the trace by a weighted trace.
\begin{prop} Let $Q\in \Cl(M, E)$ be a {\bf fixed} admissible elliptic operator with
  positive order. The exterior differential of the form 
\begin{equation} \label{eq:trQomegaj}
\omega^Q_{2j}(F) = \tr^{Q} \left(
    F(dF)^{2j}\right)
 \end{equation} 
on  ${\rm Gr}(M, E)$:
$$  d\omega^Q_{2j}= \frac{1}{{2q}} {\rm res} \left([\log Q,F] (dF)^{2k+1}
  F\right).$$ 
 is a local expression which  only  depends 
on $F$ modulo smoothing operators. 
 \end{prop}
{\bf Proof:} The locality and the dependence on $F$ modulo smoothing operators
 follow from the expression of the exterior differential in terms of a
Wodzicki residue.  To derive this expression, we mimic the finite dimensional proof, taking into account that
this time ${\rm tr}^Q$ is not cyclic:
\begin{eqnarray*}
  d\omega^Q_{2j} & = & d{\rm tr}^{Q}\left( F(dF)^{2j}\right) \\
&=& \tr^Q\left( (dF)^{2j+1}\right)\nonumber\\
&=& {\rm tr}^{Q}\left( F^2 (dF)^{2j+1}\right)\\
&=& -{\rm tr}^{Q} \left(F (dF)^{2j+1} F \right)\\
&{}& {\rm since} \quad F\, dF=-dF \, F\\
&=&\frac{1}{q} {\rm res}\left( [\log Q, F] (dF)^{2j+1}F\right) -
{\rm tr}^{Q}\left( (dF)^{2j+1} F^2\right) \\ 
&=&\frac{1}{q} {\rm res}\left( [\log Q, F] (dF)^{2j+1}F\right) -
{\rm tr}^{Q}\left( (dF)^{2j+1}\right), 
 \end{eqnarray*}
where we have used (\ref{eq:coboundary}) to write 
$$ {\rm tr}^{Q} \left([F, (dF)^{2j+1} F ] \right)= -\frac{1}{q} {\rm res} \left(F\,[ (dF)^{2j+1} F, \log
  Q]\right)=\frac{1}{q} {\rm res} \left([  F, \log
  Q]\, (dF)^{2j+1}F\right).$$
Hence 
$${\rm tr}^{Q}\left( F^2 (dF)^{2j+1}\right)= \frac{1}{2q} {\rm res}\left(
  [\log Q,
F] (dF)^{2j+1}F\right)$$ from which  the result then  follows. 
\endsquare  \\ \\ 
  Let us consider the  map \begin{eqnarray*}
\sigma:  B &\to &  \Cl_0(M, E)/ Cl_{-\infty}(M, E) \\
b&\mapsto & \bar F(b):= p\circ F(b)
\end{eqnarray*}
 where $p:\Cl_0(M, E)\to \Cl_0(M, E)/ Cl_{-\infty}(M, E)$ is the canonical projection
 map.     
In (quantum field theoretic) applications  the map
$b\mapsto \sigma(b)= \bar F_b$ can be contractible without the map $b\mapsto F_b$ being
contractible, a situation which can  occur when $B$ is contractible. To justify this, let us first observe 
that dicontinuities of $F$ give rise to jumps measured by smoothing operators
(see e.g.  \cite{Me}); indeed, since $D_b$ is a smooth family of self-adjoint elliptic
operators on a closed manifold, the discontinuities of $F_b$ are measured by
differences of 
projections $P_{b, \mu}$ over the {\it finite dimensional}
space generated by eigenvectors of $D-\lambda I$ with eigenvalues in $[0,
\mu]$ or $[\mu, 0]$ according to whether $\mu$ is positive or
negative. Finite rank projections being smoothing, it follows that the discontinuities are
measured by smoothing operators 
 so
that the projected map  $b\mapsto \bar F_b$
turns out to be 
continuous. Hence if $B$ is contractible, the map $\sigma$ is contractible.  
\begin{ex} A standard example in physics is  the case when 
\begin{itemize}
\item
the base space $B$ is a subset of connections 
in a (hermitian)  
finite-rank  vector bundle $E\to M$  over a closed (Riemannian) spin manifold
$M$. Since the space of connections in a fixed vector bundle is an affine
space, if $B$ is the whole space of connections, it is   contractible and  the map $\sigma$ is indeed contractible.
\item the operators  $D_\nabla$ (written $D_A$ with $\nabla= d+A$ in local
  coordinates)  are (twisted)  Dirac operators coupled 
to a connection $\nabla\in B$.
\end{itemize}

\end{ex}
The following theorem builds ``renormalised'' forms from the original forms $\omega^Q_{2j}$.
\begin{thm}\label{thm:renormgauge} When $\sigma(B)$ is contractible, there are  even forms $\theta^Q_{2j}$ such that 
$$ \omega^{ren}_{2j} = \omega^Q_{2j} - \theta^Q_{2j} $$ 
is closed. 
The forms $\theta^Q_{2j}$ vanish when  the order of $(dF)^{2j+1}$ is less than
-dim$\,M$. This holds in particular if  the order of $(dF)^{2j}$ is less than -dim$\,M$ in which
case  $ \omega^{ren}_{2j} = \omega^Q_{2j}= {\rm tr}(F(dF)^{2j})$ is independent of $Q$.
\end{thm}
\begin{rk}  In the case of Dirac operators parametrised by gauge connections the order of $(dF)^{2j+1}$ is less than
-dim$\,M$ for all  $j>0$ if the dimension of 
$M$ is less or equal to $-3.$  This known fact  is seen by a simple asymptotic
  expansion of the differential  $dF$. Using any fixed local trivialization of the 
underlying vector bundle $E$, we write $D_\nabla= D+A$ and $F= (D+A)
\, \vert D+A\vert^{-1}$ where $D$ is
the ordinary Dirac operator on $\R^n$ so that the infinitesimal variation $dF$
coincides up to order $1$ in $dA$ with $(D+dA)
\, \vert D+dA\vert^{-1} -D\,\vert D\vert^{-1}$. One can check that
the square of the operator $|D|^{-1} (1- \frac12 (D^{-1} dA +dA D^{-1}))$
is equal to $(D+dA)^{-2}$ up to operators of order $-3$. Hence  $dF=(D+dA) |D+dA|^{-1} = D/|D|  -\frac12 |D|^{-1} D^{-1}
\times [dA,D] +\dots $ up to operators of order $-2$.
Here one has to take into account that the commutator of $|D|$ with
$A$ is of order zero.
It follows that $D/|D|$ differs from $(D+dA)/|D+dA|$ by an operator of
order -1 so that $dF$ has order $-1$. \\
This argument fails for the case of families of metrics
 because the perturbations of Dirac operators are differential
operators of order one. It therefore does not extend to the  case of Dirac
operators parametrised by metrics since in that case the principal symbol
depends on the parameters and the differential $dF$ is a zero order
operator. 
\end{rk}
{\bf Proof  of Theorem  \ref{thm:renormgauge}:} The form  $d\omega^Q_{2j}$ being a Wodzicki residue, by Prop. 7,  it  only depends
on the projection $\bar F$  and is therefore  a pull-back  by the
  projection map $p$  of a form
  $\beta_{2j}^Q$.
The  pull-back of $\beta_{2j}^Q$  with respect to $\sigma$ 
is a closed form $\theta^Q_{2j+1}$ on $B$ which is exact since $\sigma$ is
contractible. Indeed, selecting a contraction $\sigma_t$ with $\sigma_1 =\sigma$ and $\sigma_0$ a constant map, 
we have the standard formula for the potential, $d\theta^Q_{2j} =
\theta^Q_{2j+1},$ with \begin{equation} \theta_{2j} = \frac{1}{2j+1} \int_0^1 t^{2j}\iota_{\dot\sigma_t}  \theta^Q_{2j+1}(\sigma_t) dt. \end{equation} where $\iota_X$ is the contraction by a vector field $X$ and the dot means differentiation with respect 
to the parameter $t$.\\
 When the order of $(dF)^{2j+1}$ is less than -dim$\,M$ the correction terms $\theta_{2j}^Q$ 
vanish and if  the order of $(dF)^{2j}$ is less than -dim$\,M$, the weighted trace ${\rm tr}^Q$ coincides with the usual trace so that
the naive expression $\omega^Q_{2j}$ is a closed form independent
of $Q$.
\endsquare\\ \\ 
 
The 2-form case arises in the quantum field theory gerbe \cite{CMM}. Let $B$ be a contractible  parameter 
space for Dirac operators.
 For each real number $\lambda$, let as before $U_{\lambda} \subset B$ 
be the set of parameters for which the Dirac operator $D(\lambda)=D-\lambda \,
I$ does not have $\lambda$ as an eigenvalue. 
Denote $U_{\lambda\lambda'} = U_{\lambda}\cap U_{\lambda'}$ and let $L_{\lambda\lambda'}(A)$ be the 
top exterior power of the spectral subspace $E_{\l\l'}$  defined by $\lambda <
D_A < \lambda'$ for $A\in U_{\lambda\lambda'}.$ 
The complex lines $L_{\lambda\lambda'}(A)$ form a complex line bundle $L_{\lambda\lambda'}$ over $U_{\lambda\lambda'}.$ 
For $\lambda<\lambda' <\lambda''$ we have the canonical identification
\begin{equation} L_{\lambda\lambda'} \otimes L_{\lambda'\lambda''} = L_{\lambda\lambda''} .\end{equation}
This family of line bundles defines a gerbe over $B.$ 
Since $B$ is contractible this gerbe is trivial in the sense that 
\begin{equation} L_{\lambda\lambda'} = L_{\lambda} \otimes  L_{\lambda'}^* \end{equation}
for some line bundles $L_{\lambda} \to U_{\lambda}.$ 
The curvature of $L_{\lambda\lambda'}$ is $1/2\pi$ times
\begin{equation} {\omega_2}^{\l\l'} = \tr \left( P({\l\l')} (dP(\l\l'))^2 \right)\end{equation}
where $P(\l\l')$ is the projection onto $E_{\l\l'}.$\\
 Denote as before  by $F(\l)$ the grading operator 
$(D_A-\l)/|D_A -\l|$ on $U_{\l}$ and let $P(\l)=  \frac12 (F(\l)+1)$ be the corresponding spectral projection. 
In the Hilbert-Schmidt case, when the grading operators are in 
${\rm Gr}_{\rm res}(M, E)$ one proves by a direct computation that 
\begin{equation} {\omega_2}^{\l\l'} = {\omega_2}^{\l} - {\omega_2}^{\l'} \end{equation}
on $U_{\l\l'}$ with 
\begin{equation} {\omega_2}^{\l} =  \frac18 \tr\left(F (\l)( dF(\l))^2\right)=
\tr\left(P(\l)(dP(\l))^2\right).\end{equation} 
\begin{rk} This last  equality follows from  the cyclicity of the trace on Hilbert-Schmidt operators.  
Indeed, since $d \, P(\lambda) \, P(\lambda)+ P(\lambda)\, d\,P(\lambda)=
d\, P(\lambda)$ we have  \begin{eqnarray*}\label{eq:HScase}
\tr\left( F(\l) (dF(\l))^2\right)&=& 
8 \, \tr\left( P(\l) (dP(\l))^2\right)- 4  \,  \tr\left( (dP(\l))^2\right)\\
&=&8 \, \tr\left( P(\l) (dP(\l))^2\right)- 4  \, d \tr\left( P(\l)\,
  dP(\l)\right)\\
&=&8 \, \tr\left( P(\l) (dP(\l))^2\right)+4  \, d \tr\left(dP(\l)\,
  P(\lambda)\right)\\
&=&8 \, \tr\left( P(\l) (dP(\l))^2\right)- 4  \, d \tr\left(d P(\lambda)\,
  P(\l)\right)\\
&=& 8 \, \tr\left( P(\l) (dP(\l))^2\right).
\end{eqnarray*}
\end{rk}
In the general case the forms ${\omega_2}^{\l}$ 
have to be replaced by the 'renormalised' forms as in  Theorem \ref{thm:renormgauge}. However, we still have
\begin{thm} \label{thm:linebundlecurv}
 The cocycle of forms ${\omega_2}^{\l\l'}$ is trivialized by $\frac18$ times  the forms 
\begin{equation}\label{eq:cocycle}
{\omega_2}^{\l} = \tr^{Q}\left( F(\l) (dF(\l))^2\right) - \theta^Q_2
\end{equation} 
or equivalently  by $\frac{1}{8}$ times the forms
$$ {\rho_2}^{\l} = 8 \,  \tr^{Q}\left( P(\l) (dP(\l))^2\right)  - \theta^Q_2$$
where $\theta_2^Q$ is as in  Theorem \ref{thm:renormgauge} (with $j=1$),
restricted to the open set $U_{\l}.$ In particular, 
\begin{eqnarray}\label{eq:FP}
 {\omega_2}^{\l}- {\omega_2}^{\l^\prime}
&=& {\rho_2}^{\l}-{\rho_2}^{\l^\prime}\nonumber\\
&= &
 \tr\left( F(\l) (dF(\l))^2 -F(\l^\prime) (dF(\l^\prime))^2 \right)\nonumber\\
&=&
 \tr\left( P(\l) (dP(\l))^2 -P(\l^\prime) (dP(\l^\prime))^2 \right).
\end{eqnarray}
 \end{thm}
{\bf Proof:} First observe that the correction terms $\theta^Q_2$  arising  in  the differences of the forms ${\omega_2}^{\l}$ on the intesections 
$B_{\l\l'}$ cancel: they do not depend on the parameter $\l$ since 
a change of $\l$ gives rise to finite rank perturbations of $F(\l)$ and hence
to smoothing pertubrations on which
the Wodzicki residue vanishes.\\ 
Let us show  first (\ref{eq:FP}). For $\l < \l'$ 
we have 
\begin{eqnarray*}
&{}& \tr^{Q}\left( (dP(\l))^2\right)-\tr^{Q}\left( (dP(\l^\prime))^2\right)\\
&=& 
\tr^{Q}\left( (dP(\l))^2 -(dP(\l^\prime))^2\right) \\
&=& 
\tr^{Q}\left( (dP(\l\l'))^2 + dP(\l')dP(\l\l') + dP(\l\l')dP(\l') \right) \\
&=&
\tr \left( dP(\l')dP(\l\l') + dP(\l\l') dP(\l') \right) \\
&=&
0,
\end{eqnarray*}
since  $\tr (dP)^2=0$ for any finite rank projector and
by cyclicity of the ordinary trace, from which it follows that 
$$
\tr^{Q}\left( F(\l) (dF(\l))^2\right)-\tr^{Q}\left( F(\l^\prime) (dF(\l^\prime))^2\right) = 
8 \, \left( \tr^{Q}\left( P(\l) (dP(\l))^2\right)
- \tr^{Q}\left( P(\l^\prime) (dP(\l^\prime))^2\right)\right).$$
\\ To show 
(\ref{eq:cocycle}) we 
expand  
$\omega_2^{\l}$ in powers of the
difference projection $P(\l\l')$ and  observe that the zeroth order term is equal to $\omega_2^{\l'},$ the third order term is
$\omega_2^{\l\l'}.$   The mixed terms are ordinary traces, since the operators contain the finite rank projector $P(\l\l')$ as
a factor; using the cyclicity of the trace 
 and  repeatedly $dP P = PdP =dP$ for any projector and $dP P' + P dP'=0$ for any pair of
mutually orthogonal projectors, we get 
 \begin{eqnarray*} &{}& \omega_2^{\l}  -   \omega_2^{\l'}\\
&  = & \tr \,(
 P(\l\l') dP(\l\l') dP_{\l\l'} ) + \tr\,( P(\l') dP(\l') dP(\l \l^\prime))+ \\
&+&  \tr\, (P(\l') dP(\l \l^\prime) dP(\l')) +tr\, (P(\l') dP(\l
\l^\prime)dP(\l \l^\prime)) \\
&
 +&\tr\,( P(\l \l^\prime) dP(\l') dP(\l')) +\tr\, (P(\l \l^\prime) dP(\l') dP(\l \l^\prime)) +\tr\,( P(\l \l^\prime) dP(\l \l^\prime) dP(\l')) \\
 &=&    \omega^{\l\l'}_2 - 3\tr\,\left((1-P(\l') -P(\l \l^\prime)) dP(\l')dP(\l \l^\prime)\right) .
\end{eqnarray*}
Next for any triple of mutually orthogonal projectors one has 
$\tr\,( PdP'dP'') =0,$  again by the above mentioned operator identities; applying this to $P= 1-P({\l'}) -P(\l \l^\prime), P'=P({\l'})$ and 
$P''=P(\l \l^\prime)$  we see that the mixed terms on the right-hand-side of the above equation vanish and the Theorem follows.
 \endsquare \\ 
\begin{rk}  Actually, this is just the degree 2 cohomology part of the statement that the Chern characters for direct summand in
 vector bundles add up to the Chern character of the sum; for the chosen curvature forms the statement is valid on the level of
 de Rham forms, not just for de Rham classes. 
\end{rk}
The forms $\omega_2^{\l}$ are related but not equal to the gerbe eta forms
studied in \cite{L}.  Let ${\cal G}$ be the group of smooth based group gauge
transformations acting on smooth sections $\Ci(M, E)$ as unitary operators.
On the base $B=\mathcal A/\mathcal G$ equipped with the open cover $V_{\l} = \pi(U_{\l})$ we have well-defined 
eta forms $\eta_{\l}$ such that $d\eta_{\l} = \Omega_3,$ where $\Omega_3$ is the Dixmier-Douady 3-cohomology  class 
classifying the gerbe on the base $B.$ (Here we ignore possible torsion in
cohomology and work with de Rham representatives.) Let $\pi: {\cal A}\to {\cal A}/ {\cal G}$ be the canonical projection  where
${\cal G}$  is the gauge group.  The pull-back $\pi^*(\Omega_3)$ is exact,
$d\Theta =\pi^*(\Omega_3).$ These are related to $\omega_2^{\l},$ as cohomology classes,  by
$$ [\omega_2^{\l}] = [\Theta - \pi^*(\eta_{\l})].$$
The above relation holds only in cohomology, not as a relation of forms. This is related to the fact that the 
difference of the pull-back forms $\pi^*(\eta_{\l})$ must vanish in  gauge directions whereas the difference 
$\omega^{\l}_2 -\omega^{\l'}_2= \omega^{\l\l'}_2$ is nonvanishing even in gauge directions. However, we have 

\begin{prop} Restricted to gauge orbits, the forms $\omega^{\l\l'}_2$ vanish as cohomology classes. More precisely,
$\omega^{\l\l'}_2 = d\omega^{\l\l'}_1$ on gauge orbits, where
$$\omega_1^{\l\l'}(X) = -\tr\left(\, X P(\l \l^\prime)\right).$$ 
Here we can use ordinary trace since $P(\l \l^\prime)$ has finite rank; $X$ is an  element  of the Lie algebra of $\mathcal G$ 
acting as multiplication operator in the Hilbert space $H.$  \end{prop} 

{\bf Proof:} The gauge group acts on projections by conjugation $P\mapsto gPg^{-1}$ and so 
\begin{eqnarray*}&{}& (d\omega^{\l\l'}_1)(X,Y)\\
 &=& \mathcal L_X \omega^{\l\l'}_1(Y) - \mathcal L_Y \omega^{\l\l'}_1(X) - \omega^{\l\l'}_1([X,Y]) \\
&=& -\tr\left(Y [P(\l \l^\prime),X]\right) +\tr\left( X [P(\l \l^\prime),Y] \right)+\tr\left( [X,Y]P(\l \l^\prime)\right) \\
&=& -\tr\left( [X,Y]P(\l \l^\prime)\right)\\
&=& \tr\left(P(\l \l^\prime) [[P(\l \l^\prime), X],[P(\l \l^\prime) ,Y]]\right)\\
&=& \omega^{\l\l'}_2(X,Y), \end{eqnarray*} 
where in the last step we have used the projection property $P^2=P$ and the cyclicity of trace 
for finite rank operators; $\mathcal L_X$ is the Lie derivative by vector field along gauge orbits corresponding to 
the conjugation action of the group $\mathcal G.$  \endsquare

\begin{rk} \rm A similar modification can be made for the gauge action on the local forms $\omega^{\l}_2$ to show 
that the action is consistent on overlaps. Again, restricting to gauge orbits, using $F^2=1$ and rearranging terms, one can write 
\begin{eqnarray*}&{}& \tr^{Q}\left( F(\l)[[F(\l), X], [F(\l),Y] \right)\\
&=& -4\tr^{Q} \left([X,Y] F(\l) \right)
+ 2\tr^{Q}\left(2[XF(\l) ,Y]+[F(\l)X F(\l)Y,F(\l)] 
+[XY,F(\l)]\right) \\ 
&=& -4\tr^{Q}\left([X,Y]F(\l)\right) + 2\res \left( 2[\log Q, XF(\l)]Y+ [\log
  Q, XY]F(\l)+[\log Q, F(\l)XF(\l)Y]F(\l)  \right)\end{eqnarray*}
where on the right the first term is a trivial cocycle and the rest, being a residue, does not
depend on finite rank perturbations and in particular not on the parameter $\l.$ 

Here we have only discussed the cocycles of degree $2$ because they are the most relevant in gauge theory; it is clear that 
similar computations can be performed with the higher cocycles. 
\\ The 2-forms $\omega_2^{\l}$ are directly 'seen' in quantum field theory in the following way. These forms appear as curvature forms
of local vacuum line bundles for fermion field in gauge background, \cite{CMM, EM}. The gauge action on gauge connections
lifts to an action of an extension of the group of gauge transformations  on the local line bundles. On the Lie algebra level, the 2-cocycle describing the
Lie algebra extension(``hamiltonian anomaly'' \cite{Mi})  is just the curvature form evaluated in the gauge directions. In the case of fields in one space dimension,
this extension (for a simple compact gauge group) defines an affine Kac-Moody algebra.  
\end{rk}

\bibliographystyle{plain}

\end{document}